\documentclass{amsart}
\usepackage{hyperref}
\usepackage{amssymb}
\usepackage{amsmath}
\usepackage{amsthm}
\usepackage{amsfonts}
\usepackage{graphicx}
\usepackage{amscd}
\usepackage{mathrsfs}
\usepackage{euscript}

\allowdisplaybreaks
\makeatletter
\@ifundefined{pdfoutput}
{
\DeclareGraphicsRule{.wmf}{bmp}{}{}
\DeclareGraphicsRule{.jpg}{bmp}{}{}
\DeclareGraphicsRule{.png}{bmp}{}{}
\DeclareGraphicsRule{.cdr}{bmp}{}{}
\DeclareGraphicsRule{.gif}{bmp}{}{}
}
{
\ifnum\pdfoutput=0\relax
\DeclareGraphicsRule{.wmf}{bmp}{}{}
\DeclareGraphicsRule{.jpg}{bmp}{}{}
\DeclareGraphicsRule{.png}{bmp}{}{}
\DeclareGraphicsRule{.cdr}{bmp}{}{}
\DeclareGraphicsRule{.gif}{bmp}{}{}
\fi
\ifnum\pdfoutput=1\relax
\DeclareGraphicsRule{.eps}{pdf}{}{}
\DeclareGraphicsRule{.wmf}{jpg}{}{}
\fi
}
\makeatother
\newtheorem{thm}{Theorem}[section]

\newtheorem{lem}[thm]{Lemma}
\newtheorem{nota}[thm]{Notation}
\newtheorem{prop}[thm]{Proposition}

\theoremstyle{definition}

\newtheorem{defn}[thm]{Definition}

\newtheorem{example}{Example}[section]
\newtheorem{remark}[thm]{Remark}
\newtheorem{notation}[thm]{Notation}
\numberwithin{equation}{section}

\title[Stochastic integrals and BM on abstract nilpotent Lie groups]
{Stochastic integrals and Brownian motion on abstract nilpotent Lie groups}
\author[Melcher]{Tai Melcher{$^*$}}
\thanks{\footnotemark {$^*$} This research was supported in part by NSF
Grants DMS-0907293 and DMS-1255574.}
\address{Department of Mathematics\\
University of Virginia\\ Charlottesville, VA 22903 USA}
\email{melcher@virginia.edu}


\begin{document}

\begin{abstract}
We construct a class of iterated stochastic integrals with respect to Brownian motion on an abstract Wiener space which allows for the definition of Brownian motions on a general class of infinite-dimensional nilpotent Lie groups based on abstract Wiener spaces. We then prove that a Cameron--Martin type quasi-invariance result holds for the associated heat kernel measures in the non-degenerate case, and give estimates on the associated Radon--Nikodym derivative.  We also prove that a log Sobolev estimate holds in this setting.
\end{abstract}

\maketitle
\tableofcontents

\section{Introduction\label{s.1}}


The construction of diffusions on infinite-dimensional manifolds and the study of the regularity properties of their induced measures have been a topic of great interest for at least the past 50 years; see for example \cite{DalettskiSnaiderman1969,Kuo1971,Kuo1972,Elworthy1975,FdlP1995,ADK2003,AM2006,Malliavin08}, although many other references exist.  The purpose of the present paper is to construct diffusions on a general class of infinite-dimensional nilpotent Lie groups, and to show that the associated heat kernel measures are quasi-invariant under  appropriate translations.  We demonstrate that this class of groups is quite rich. We focus here on the elliptic setting, but comment that, as nilpotent groups are standard first models for studying hypoellipticity, examples of infinite-dimensional versions of such spaces are important for the more general study of hypoellipticity in infinite dimensions.  This is an area of great interest, and is of particular relevance in the study of stochastic PDEs and their applications \cite{AiraultMalliavin2004,BakhtinMattingly2007,BaudoinTeichmann2005,HairerMattingly2011,Malliavin1990,MattinglyPardoux2006}. The present paper studies heat kernel measures for elliptic diffusions on these spaces, which is a necessary precursor to understanding the degenerate case.

\subsection{Main results}
Let $(\mathfrak{g},\mathfrak{g}_{CM},\mu)$ denote an abstract Wiener
space, where $\mathfrak{g}$ is a Banach space equipped with a
centered non-degenerate Gaussian measure $\mu$ and associated
Cameron--Martin Hilbert space $\mathfrak{g}_{CM}$.  We will assume
that $\mathfrak{g}_{CM}$ additionally carries a nilpotent Lie algebra
structure $[\cdot,\cdot]: \mathfrak{g}_{CM}\times \mathfrak{g}_{CM}
\rightarrow \mathfrak{g}_{CM}$, and we will further assume that
this Lie bracket is Hilbert-Schmidt.  We will call such a space an {\it
abstract nilpotent Lie algebra}. Via the standard
Baker-Campbell-Hausdorff-Dynkin formula, we may then equip
$\mathfrak{g}_{CM}$ with an explicit group operation under which
$\mathfrak{g}_{CM}$ becomes an infinite-dimensional group.  When
thought of as a group, we will denote this space by $G_{CM}$.  
We equip $G_{CM}$ with the left-invariant Riemannian metric which
agrees with the inner product on $\mathfrak{g}_{CM}\cong T_eG_{CM}$, and we denote the
Riemannian distance by $d$.  
Note that, despite the
use of the notation $\mathfrak{g}$, we do not assume that the Lie bracket
structure extends to $\mathfrak{g}$, and so this space is not necessarily a
Lie algebra or a Lie group.  Still, when $\mathfrak{g}$ is playing a role
typically played by the group, we will denote $\mathfrak{g}$ by $G$.

If $\mathfrak{g}$ were a Lie algebra (and thus carried an associated group
structure), we could construct a Brownian motion on
$G$ as the solution to the Stratonovich stochastic differential equation
\[ \delta g_t = g_t\delta B_t := L_{g_t*}\delta B_t, \text{ with }
	g_0=\mathbf{e}=(0,0), \]
where $L_x$ is left translation by $x\in G$ and $\{B_t\}_{t\ge0}$ is a
standard Brownian motion on $\mathfrak{g}$ (as a Banach space) with $\mathrm{Law}(B_1)=\mu$.
In finite dimensions, the solution to this stochastic differential equation may be obtained
explicitly as a formula involving the Lie bracket.  In particular, for $t>0$
and $n\in\mathbb{N}$, let $\Delta_n(t)$ denote the simplex in
$\mathbb{R}^n$ given by
\[ \{s=(s_1,\cdots,s_n)\in\mathbb{R}^n: 0<s_1<s_2<\cdots<s_n<t\}. \]
Let $\mathcal{S}_n$ denote the permutation group on $(1,\cdots,n)$,
and, for each $\sigma\in \mathcal{S}_n$, let $e(\sigma)$ denote the
number of ``errors'' in the ordering
$(\sigma(1),\sigma(2),\cdots,\sigma(n))$, that is, $e(\sigma)=\#
\{j<n: \sigma(j)>\sigma(j+1)\}$.
Then the Brownian motion on $G$ could be written as 
\begin{equation}
\label{e.ibm}
g_t = \sum_{n=1}^{r-1} \sum_{\sigma\in\mathcal{S}_n}  
	\left( (-1)^{e(\sigma)}\bigg/ n^2 
	\begin{bmatrix} n-1 \\ e(\sigma) \end{bmatrix}\right) 
	\int_{\Delta_n(t)} 
	[ [\cdots[\delta B_{s_{\sigma(1)}},\delta B_{s_{\sigma(2)}}],\cdots], 
	\delta B_{s_{\sigma(n)}}],
\end{equation}
where this sum is finite under the assumed nilpotence.  
An obstacle to the development of a general theory of stochastic differential equations on infinite-dimensional Banach spaces is the lack of smoothness of the norm in a general Banach space which is necessary to define a stochastic integral on it.
Still, in Section \ref{s.mult}, we prove a general result to define a class of iterated stochastic integrals with respect to Brownian motion on the Banach space $\frak{g}$ that includes the expression above.  Additionally, we show that one may make
sense of the above expression when the Lie bracket on $\mathfrak{g}_{CM}$ does not necessarily extend to $\mathfrak{g}$. Thus we are able to define a ``group Brownian motion'' $\{g_t\}_{t\ge0}$ on $G$ via (\ref{e.ibm}). We let $\nu_t:=\mathrm{Law}(g_t)$ be the heat kernel measure on $G$.

In particular, the integrals above are defined as a limit of stochastic integrals on finite-dimensional subgroups $G_\pi$ of $G_{CM}$. We show that these $G_\pi$ are nice  in the sense that they approximate $G_{CM}$ and that there exists a uniform lower bound on their Ricci curvatures.  

Using these results, we are able to prove the following main theorem.

\begin{thm}
\label{t.quasi}
For $h\in G_{CM}$, let $L_h,R_h:G_{CM}\rightarrow G_{CM}$ denote left and
right translation by $h$, respectively. Then $L_h$ and $R_h$ define measurable
transformations on $G$, and for all $T>0$, $\nu_t\circ L_h^{-1}$ and $\nu_t\circ
R_h^{-1}$ are absolutely continuous with respect to $\nu_t$.  Let 
\[ J_t^l(h,\cdot) := \frac{d(\nu_t\circ L_h^{-1})}{d\nu_t} \qquad \text{ and }
	\qquad J_t^r(h,\cdot) := \frac{d(\nu_t\circ R_h^{-1})}{d\nu_t} \]
be the Radon-Nikodym derivatives, $k$ be the uniform lower bound on the Ricci
curvatures of the finite-dimensional approximation groups $G_\pi$ 
and 
\[ c(t) := \frac{t}{e^t-1}, \qquad \text{ for all } t\in\mathbb{R}, \]
with the convention that $c(0)=1$.  Then, for all $p\in[1,\infty)$,
$J_t^l(h,\cdot),J_t^r(h,\cdot)\in L^p(\nu_t)$ and both satisfy the estimate
\[ \|J_t^*(h,\cdot)\|_{L^p(\nu_t)} \le \exp \left( \frac{c(kt)(p-1)}{2t}
	d(\mathbf{e},h)^2\right), \]
where $*=l$ or $*=r$.
\end{thm}

The fact that one may define a measurable left or right action on $G$ by an
element of $G_{CM}$ is discussed in Section \ref{s.mga}.  The lower bound on
the Ricci curvature is proved in Proposition \ref{p.Ric}.

\subsection{Discussion}

The present paper builds on the previous work in \cite{DriverGordina2008} and
\cite{Melcher2009}, significantly generalizing
these previous works in several ways.
In particular, the paper \cite{Melcher2009} considered
analogous results for ``semi-infinite Lie groups'', which are 
infinite-dimensional nilpotent Lie groups
constructed as extensions of finite-dimensional nilpotent Lie groups
$\mathfrak{v}$ by an
infinite-dimensional abstract Wiener space (see Example \ref{ex.ext}).  
At several points in the analysis
there, the fact that $\mathrm{dim}(\mathfrak{v})<\infty$ was used in a
critical way.  In particular, it was used to show that the stochastic
integrals defining the Brownian motion on $G$ as in equation
(\ref{e.ibm}) were well-defined.  
In the present paper, we have removed this restriction, as well as removing the
``stratified'' structure implicit in the construction as a Lie group extension.

Again, we note that, despite the
use of the notation $\mathfrak{g}$, it is not assumed that the Lie bracket
structure on $\mathfrak{g}_{CM}$ extends to $\mathfrak{g}$, and so $\mathfrak{g}$ itself is not necessarily a
Lie algebra or Lie group.  In \cite{DriverGordina2008} and \cite{Melcher2009}, it
was assumed that the Lie bracket was a continuous map defined on
$\mathfrak{g}$.  However, it turns out that
the group construction on $\mathfrak{g}_{CM}$ is  the only
necessary structure for the subsequent analysis.
As is usual for the infinite-dimensional setting,
while the heat kernel measure is itself supported on the larger space
$\mathfrak{g}$, its critical analysis depends more on the structure of $\mathfrak{g}_{CM}$.
Still, as was originally done in \cite{DriverGordina2008} and then in \cite{Melcher2009},
one may instead define an abstract nilpotent Lie algebra starting 
with a continuous nilpotent bracket
$[\cdot,\cdot]:\mathfrak{g}\times\mathfrak{g}\rightarrow\mathfrak{g}$.  For
example, in the event that $\mathfrak{g}=W\times\mathfrak{v}$ where
$\mathfrak{v}$ is a finite-dimensional Lie algebra and
$[\cdot,\cdot]:\mathfrak{g}\times\mathfrak{g}\rightarrow\mathfrak{v}$, it is
well-known that this implies that the restriction of the bracket to
$\mathfrak{g}_{CM}=H\times\mathfrak{v}$ is Hilbert-Schmidt.  (For
any continuous bilinear $\omega:W\times W\rightarrow K$ where $K$ is a Hilbert
space, one has that $\|\omega\|_{(H^{\otimes 2})^*\otimes K}<\infty$; this
follows for example from Corollary 4.4 of \cite{Kuo75}.)
More generally, in order for the subsequent theory to make sense, one would
naturally need
to require that $\mathfrak{g}_{CM}$ be a Lie subalgebra of $\mathfrak{g}$, that is, for
the restriction of the Lie bracket to $\mathfrak{g}_{CM}$ to preserve
$\mathfrak{g}_{CM}$.  As the proofs in the sequel rely strongly on the bracket
being Hilbert-Schmidt, it would be then necessary to add the Hilbert-Schmidt hypothesis as it 
does not follow immediately if one only assumes a continuous bracket on
$\mathfrak{g}$ which preserves $\mathfrak{g}_{CM}$.

Additionally, the spaces studied in the present paper are well-designed for the study of
infinite-dimensional hypoelliptic heat kernel measures, and there has already
been progress on proving quasi-invariance and stronger smoothness properties
for these measures in the simplest case of a step two Lie algebra with
finite-dimensional center; see \cite{BGM2013} and
\cite{DEM2013}. More generally, the paper \cite{Pickrell2011} explores related interesting lines of inquiry for heat kernel measures on infinite-dimensional groups, largely in the context of groups of maps from manifolds to Lie groups.

{\it Acknowledgements. } The author thanks Bruce Driver and Nathaniel Eldredge for
helpful conversations during the writing of this paper.

\section{Iterated It\^o integrals}
\label{s.mult}

Recall the standard construction of abstract Wiener spaces.  
Suppose that $W$ is a real separable Banach space and $\mathcal{B}_{W}$ is
the Borel $\sigma$-algebra on $W$.

\begin{defn}
\label{d.2.1} 
A measure $\mu$ on $(W,\mathcal{B}_{W})$ is called a (mean zero,
non-degenerate) {\em Gaussian measure} provided that its characteristic
functional is given by
\[
\hat{\mu}(u) := \int_W e^{iu(x)} d\mu(x)
	= e^{-\frac{1}{2}q(u,u)}, \qquad \text{ for all } u\in W^*,
\]
for $q=q_\mu:W^*\times W^*\rightarrow\mathbb{R}$ a symmetric, positive
definite quadratic form.
That is, $q$ is a real inner product on $W^*$.
\end{defn}

\begin{thm}
\label{t.2.3}
Let $\mu$ be a Gaussian measure on $W$.  
For $w\in W$, let
\[ 
\|w\|_H := \sup\limits_{u\in W^*\setminus\{0\}}\frac{|u(w)|}{\sqrt{q(u,u)}},
\]
and define the {\em Cameron--Martin subspace} $H\subset W$ by
\[ H := \{h\in W : \|h\|_H < \infty\}. \]
Then $H$ is a dense subspace of $W$, and there exists a unique inner
product $\langle\cdot,\cdot\rangle_H$ on $H$ such that $\|h\|_H^2
= \langle h,h\rangle_H$ for all $h\in H$, and $H$ is a separable
Hilbert space with respect to this inner product.  For any $h\in
H$, $\|h\|_W \le C \|h\|_H$ for some $C<\infty$.  
\end{thm}

Alternatively, given $W$ a
real separable Banach space and $H$ a real separable Hilbert space 
continuously embedded in $W$ as a dense subspace, then for each $w^*\in W^*$ there
exists a unique $h_{w^*}\in H$ such that $\langle h,w^*\rangle = \langle h,
h_{w^*}\rangle_H$ for all $h\in H$.  Then $W^*\ni w^*\mapsto h_{w^*}\in
H$ is continuous, linear, and one-to-one with a dense range
\begin{equation}
\label{e.H*} 
H_* := \{h_{w^*}:w^*\in W^*\},
\end{equation}
and $W^*\ni w^*\mapsto h_{w^*}\in
W$ is continuous. A Gaussian measure on $W$ is a Borel
probability measure $\mu$ such that, for each $w^*\in W^*$, the random
variable $w\mapsto\langle w,w^*\rangle$ under $\mu$ is a centered Gaussian
with variance $\|h_{w^*}\|_H^2$.

Suppose that $P:H\rightarrow H$ is a finite rank orthogonal projection
such that $PH\subset H_*$. Let $\{h_j\}_{j=1}^m$ be an orthonormal basis for
$PH$.  Then we may extend $P$ to a (unique) continuous operator
from $W\rightarrow H$ (still denoted by $P$) by letting
\begin{equation}
\label{e.proj}
Pw := \sum_{j=1}^m \langle w, h_j\rangle_H  h_j
\end{equation}
for all $w\in W$.  
\begin{nota}
\label{n.proj}
Let $\mathrm{Proj}(W)$ denote the collection of finite rank projections
on $W$ such that $PW\subset H_*$ and $P|_H:H\rightarrow H$ is an orthogonal
projection, that is, $P$ has the form given in equation
\eqref{e.proj}.
\end{nota}

Let $\{B_t\}_{t\ge0}$ be a Brownian motion on $W$ with variance determined by
\[
\mathbb{E}\left[\langle B_s,h\rangle
		_{H} 
		\langle B_t,k\rangle
		_H\right] 
	= \langle h,k \rangle_{H} \min(s,t),
\]
for all $s,t\ge0$ and $h,k\in H_*$, where $H_*$ is as in (\ref{e.H*}).  Note that for any $P\in\mathrm{Proj}(W)$, $PB$ is a Brownian motion on $PH$. In the rest of this section, we will verify the existence of martingales defined as certain iterated stochastic integrals with respect to $B_t$.

The following is Proposition 4.1 of \cite{Melcher2009}.  Note that again this
was stated in the context where $H=\mathfrak{g}_{CM}$ was a ``semi-infinite Lie algebra'', but 
a brief inspection of the proof shows that this is a general statement about stochastic
integrals on Hilbert spaces.

\begin{prop}
\label{p.int1}
Let $\{P_m\}_{m=1}^\infty\subset\mathrm{Proj}(W)$ such that $P_m|_H\uparrow
I_H$.  Then, for $\xi\in L^2(\Delta_n(t),H^{\otimes n})$ a continuous 
mapping, let
\begin{align*}
J_n^m(\xi)_t &:= \int_{\Delta_n(t)} \langle P_m^{\otimes n} \xi(s), dB_{s_1}
		\otimes\cdots\otimes dB_{s_n}
		\rangle_{H^{\otimes n}} \\
	&= \int_{\Delta_n(t)} \langle \xi(s), dP_m B_{s_1}
	\otimes\cdots\otimes dP_m B_{s_n}
	\rangle_{H^{\otimes n}}. 
\end{align*}
Then $\{J_n^m(\xi)_t\}_{t\ge0}$ is a continuous $L^2$-martingale, 
and there exists a
continuous $L^2$-martingale $\{J_n(\xi)_t\}_{t\ge0}$ such that
\begin{equation} 
\label{e.Jnm}
\lim_{m\rightarrow\infty} \mathbb{E}\left[ \sup_{\tau\le t} 
	|J_n^m(\xi)_\tau-J_n(\xi)_\tau|^2 \right] = 0
\end{equation}
and  
\begin{equation} 
\label{e.xi}
\mathbb{E}|J_n(\xi)_t|^2\le
	\|\xi\|^2_{L^2(\Delta_n(t),H^{\otimes n})}
\end{equation}
for all $t<\infty$.  The process
$J_n(\xi)$ is well-defined independent of the choice of increasing
orthogonal projections $\{P_m\}_{m=1}^\infty$ into $H_*$, and so will be denoted by
\[ J_n(\xi)_t 
	= \int_{\Delta_n(t)} \langle \xi(s), dB_{s_1}\otimes\cdots\otimes dB_{s_n}
	\rangle_{H^{\otimes n}}.
\]
\end{prop}

Now we may use this result to define stochastic integrals taking values in
another Hilbert space $K$.

\begin{prop}
\label{p.int2}
Let $K$ be a Hilbert space and $F\in L^2(\Delta_n(t),(H^{\otimes
n})^*\otimes K)$ be a continuous map.  That is, $F:\Delta_n(t)\times H^{\otimes
n}\rightarrow K$ is a map continuous in $s$ and linear on $H^{\otimes n}$ such that
\[ \int_{\Delta_n(t)} \|F(s)\|_2^2\,ds 
	= \int_{\Delta_n(t)} \sum_{j_1,\ldots,j_n=1}^\infty
\|F(s)(h_{j_1}\otimes\cdots\otimes h_{j_n})\|_K^2\,ds <\infty. \]
Then
\[ J_n^m(F)_t := \int_{\Delta_n(t)} F(dP_m B_{s_1}
	\otimes\cdots\otimes dP_m B_{s_n}) 
\]
is a continuous $K$-valued $L^2$-martingale, and there exists a continuous $K$-valued $L^2$-martingale
$\{J_n(F)_t\}_{t\ge0}$ such that
\begin{equation}
\label{e.intF}
\lim_{m\rightarrow\infty} \mathbb{E}\left[ \sup_{\tau\le t} 
	\|J_n^m(F)_\tau-J_n(F)_\tau\|_{K}^2 \right] = 0,
\end{equation}
for all $t<\infty$.  The martingale $J_n(F)$ is well-defined independent of 
the choice of orthogonal projections, 
and thus will be denoted by
\[
J_n(F)_t = \int_{\Delta_n(t)} F(dB_{s_1}\otimes\cdots\otimes dB_{s_n}).
\]
\end{prop}

\begin{proof}
Let $\{e_j\}_{j=1}^{\infty}$ be an orthonormal basis of $K$. 
Since $\langle F(s)(\cdot), e_j \rangle$ is linear on
$H^{\otimes n}$, for each $s$ there exists
$\xi_j(s)\in H^{\otimes n}$ such that
\begin{equation}
\label{e.xij} 
\langle \xi_j(s), k_1\otimes\cdots\otimes k_n \rangle 
	= \langle F(s)(k_1\otimes\cdots\otimes k_n), e_j \rangle. 
\end{equation}
If $\xi_j:\Delta_n(t)\rightarrow H^{\otimes n}$ 
is defined by equation (\ref{e.xij}), 
then clearly $\xi_j\in L^2(\Delta_n(t),H^{\otimes n})$ and in particular
\begin{align*}
\|F\|_{L^2(\Delta_n(t)\times H^{\otimes n},K)}^2
	&=  \sum_{j=1}^\infty \|\xi_j\|_{L^2(\Delta_n(t),H^{\otimes n})} <\infty.
\end{align*}
Thus, for $J_n(\xi_j)$ as defined in Proposition \ref{p.int1},
\begin{align*}
\mathbb{E}\left[ \sum_{j=1}^\infty |J_n(\xi_j)_t|^2\right] 
	&\le \frac{t^n}{n!}\mathbb{E}\left[ \int_{\Delta_n(t)} \sum_{j=1}^\infty |\langle \xi_j(s),
		dB_{s_1}\otimes\cdots\otimes dB_{s_n}\rangle_{H^{\otimes n}}|^2
		\right] \\
	&= \frac{t^n}{n!}\sum_{j=1}^\infty \|\xi_j\|^2_{L^2(\Delta_n(t),H^{\otimes n})}
	<\infty,
\end{align*}
and so we may write
\begin{align*} 
\sum_{j=1}^\infty J_n(\xi_j)_t e_j 
	&= \sum_{j=1}^\infty \int_{\Delta_n(t)} 
		\langle \xi_j(s),
		dB_{s_1}\otimes\cdots\otimes dB_{s_n}\rangle_{H^{\otimes n}} e_j  \\
	&= \int_{\Delta_n(t)} \sum_{j=1}^\infty
		\langle F(s)(
		dB_{s_1}\otimes\cdots\otimes dB_{s_n}),e_j\rangle_K e_j \\
	&= \int_{\Delta_n(t)} 
		F(s)(dB_{s_1}\otimes\cdots\otimes dB_{s_n}).
\end{align*}
Thus, taking $J_n(F)_t := \sum_{j=1}^\infty J_n(\xi_j)_t e_j$,
we also have that 
\begin{align*}
\mathbb{E}\|J_n(F)_t - J_n^m(F)_t\|_K^2
	&= \mathbb{E}\left[\sum_{j=1}^\infty
		|J_n(\xi_j)_t-J_n^m(\xi_j)_t|^2\right]
	\rightarrow 0
\end{align*}
as $m\rightarrow\infty$ by (\ref{e.Jnm}) and dominated convergence since
\[ \mathbb{E}|J_n(\xi_j)_t-J_n^m(\xi_j)_t|^2 
	\le 4\|\xi_j\|_{L^2(\Delta_n(t),H^{\otimes n})}^2 \]
by (\ref{e.xi}).  Then equation (\ref{e.intF}) holds by Doob's maximal
inequality.
\end{proof}

Note that the preceding results then imply that one may define the above
stochastic integrals with respect to {\it any} increasing sequence of
orthogonal projections -- that is, we need not require that the projections
extend continuously to $W$.


\begin{prop}
\label{p.bad}
Let $V$ be an arbitrary finite-dimensional subspace of $H$, and let
$\pi:H\rightarrow V$ denote orthogonal projection onto $V$.  Then for any
Hilbert space $K$ and  $F\in L^2(\Delta_n(t),(H^{\otimes
n})^*\otimes K)$ a continuous map, the
stochastic integral 
\[ J_n^\pi(F)_t := \int_{\Delta_n(t)} F(d\pi B_{s_1}
	\otimes\cdots\otimes d\pi B_{s_n}) 
\]
is well-defined, and $\{J_n^\pi(F)_t\}_{t\ge0}$
is a continuous $K$-valued $L^2$-martingale.  Moreover, if $V_m$ is an increasing sequence of
finite-dimensional subspaces of $H$
such that the corresponding orthogonal projections $\pi_m\uparrow I_H$, then 
\begin{equation*} 
\lim_{m\rightarrow\infty} \mathbb{E}\left[ \sup_{\tau\le t} 
	\|J_n^{\pi_m}(F)_\tau-J_n(F)_\tau\|^2 \right] = 0,
\end{equation*}
where $J_n(F)$ is as defined in Proposition \ref{p.int2}.
\end{prop}

\begin{proof}
First consider the case that $K=\mathbb{R}$, and thus $F(s)=\langle
\xi(s),\cdot\rangle$ for a continuous $\xi\in L^2(\Delta_n(t),H^{\otimes n})$.
Since $\pi^{\otimes n}\xi\in L^2(\Delta_n(t),H^{\otimes n})$, the definition
of $J_n^\pi(\xi)=J_n(\pi^{\otimes n}\xi)$ follows from Propostion \ref{p.int1}.  Moreover, by equation
(\ref{e.xi}), 
\begin{align*}
\mathbb{E}|J_n^{\pi_m}(\xi)_t - J_n(\xi)_t|^2
	&= \mathbb{E}|J_n(\pi_m^{\otimes n}\xi)_t - J_n(\xi)_t|^2 
	= \mathbb{E}|J_n(\pi_m^{\otimes n}\xi-\xi)_t|^2 \\
	&\le \|\pi_m^{\otimes n}\xi-\xi\|_{L^2(\Delta_n(t),H^{\otimes n})}
	\rightarrow 0
\end{align*}
as $m\rightarrow\infty$.  Now the proof for general $F$ follows just as in
Proposition \ref{p.int2}.
\end{proof}

\section{Abstract nilpotent Lie algebras and groups}
\label{s.prelim}

\begin{defn}
\label{d.semi}
Let $(\mathfrak{g},\mathfrak{g}_{CM},\mu)$ be an abstract
Wiener space such that $\mathfrak{g}_{CM}$ is equipped with a nilpotent
Hilbert-Schmidt Lie bracket.  Then we will call
$(\mathfrak{g},\mathfrak{g}_{CM},\mu)$ an {\it abstract nilpotent Lie
algebra}.
\end{defn}

The Baker-Campbell-Hausdorff-Dynkin formula implies that
\[ \log(e^A e^B) = A+B+\sum_{k=1}^{r-1} 
		\sum_{(n,m)\in\mathcal{I}_k}  
		a_{n,m}^k\mathrm{ad}_A^{n_1} \mathrm{ad}_B^{m_1} \cdots
		\mathrm{ad}_A^{n_k} \mathrm{ad}_B^{m_k} A,
\]
for all $A,B\in\mathfrak{g}_{CM}$, where 
\begin{equation*} 
a_{n,m}^k := \frac{(-1)^k}{(k+1)m!n!(|n|+1)}, 
\end{equation*}
$\mathcal{I}_k := \{(n,m)\in\mathbb{Z}_+^k\times\mathbb{Z}_+^k : 
n_i+m_i>0 \text{ for all } 1\le i\le k \}$, and for each multi-index
$n\in\mathbb{Z}_+^k$,
\[ n!= n_1!\cdots n_k! \quad \text{ and } \quad |n|=n_1+\cdots+n_k; \]
see, for example, \cite{duiskolk}.  If $\mathfrak{g}_{CM}$ is nilpotent of step
$r$, then
\[ \mathrm{ad}_A^{n_1} \mathrm{ad}_B^{m_1} \cdots
		\mathrm{ad}_A^{n_k} \mathrm{ad}_B^{m_k} A = 0 \quad
\text{if } |n|+|m|\ge r. \]
for $A,B\in\mathfrak{g}_{CM}$.  
Since $\mathfrak{g}_{CM}$ is simply connected and nilpotent, 
the exponential map is a global diffeomorphism (see, for
example, Theorems 3.6.2 of \cite{Varadarajan} or 1.2.1 of \cite{CorGrn90}). 
In particular, we may view $\mathfrak{g}_{CM}$ as both a Lie algebra and Lie group, and 
one may verify that
\begin{align}
\label{e.mult}
g\cdot h &= g+h+\sum_{k=1}^{r-1} 
		\sum_{(n,m)\in\mathcal{I}_k}  
		a_{n,m}^k\mathrm{ad}_g^{n_1} \mathrm{ad}_h^{m_1} \cdots
		\mathrm{ad}_g^{n_k} \mathrm{ad}_h^{m_k} g
\end{align}
defines a group structure on $\mathfrak{g}_{CM}$.  Note that $g^{-1}=-g$ and 
the identity $\mathbf{e}=(0,0)$.

\begin{defn}
When we wish to emphasize the group structure on $\mathfrak{g}_{CM}$, we will
denote $\mathfrak{g}_{CM}$ by $G_{CM}$.
\end{defn}

\begin{lem}
\label{l.cts}
The Banach space topology on $\mathfrak{g}_{CM}$ makes $G_{CM}$ into a topological group.
\end{lem}
\begin{proof}
Since $\mathfrak{g}_{CM}$ is a topological vector space, 
$g\mapsto g^{-1}=-g$ and $(g_1,g_2)\mapsto g_1+g_2$ are continuous by
definition.  The map $(g_1,g_2)\mapsto [g_1,g_2]$ is continuous in 
the $\mathfrak{g}_{CM}$ topology by the boundedness of the Lie bracket.
It then follows from (\ref{e.mult}) that $(g_1,g_2)\mapsto g_1\cdot g_2$
is continuous as well.
\end{proof}

\label{s.mga}\subsection{Measurable group actions on $G$}

As discussed in the introduction, given a Hilbert-Schmidt Lie bracket on
$\mathfrak{g}_{CM}$ and a subsequently defined group operation on $G_{CM}$,
one may define a measurable action on $G$ by left or right multiplication by
an element of $G_{CM}$.  

In particular, let $\{e_n\}_{n=1}^\infty$ be an orthonormal basis of $\mathfrak{g}_{CM}$.  
For now, fix $n$ and consider the mapping $\mathfrak{g}_{CM}\rightarrow \mathfrak{g}_{CM}^*$ given by $h\mapsto \langle
\mathrm{ad}_h\cdot,e_n\rangle$.  Then this is a continuous linear map on $\mathfrak{g}_{CM}$
and in the usual way we may make the identification of
$\mathfrak{g}_{CM}^*\cong \mathfrak{g}_{CM}$ so that we
define the operator $A_n:\mathfrak{g}_{CM}\rightarrow \mathfrak{g}_{CM}$ given by
\[ \langle A_nh, k\rangle = \langle \mathrm{ad}_h k,e_n\rangle; \]
in particular, $A_nh=\mathrm{ad}_h^*e_n$.  Note that, for any $h,k\in \mathfrak{g}_{CM}$
\[ \langle A_n^*h,k\rangle  
	= \langle A_nk,h\rangle 
	= \langle \mathrm{ad}_k^* e_n,h\rangle
	= \langle e_n , \mathrm{ad}_k h \rangle
	= - \langle e_n , \mathrm{ad}_h k \rangle
	= - \langle \mathrm{ad}_h^* e_n,k\rangle \]
and thus $A_n^*=-A_n$.
Now  fix $h\in G_{CM}=\mathfrak{g}_{CM}$.  Then for $\mathrm{ad}_h:\mathfrak{g}_{CM}\rightarrow \mathfrak{g}_{CM}$ we may write
\[
\mathrm{ad}_h k = \sum_n \langle \mathrm{ad}_hk, e_n\rangle e_n
	= \sum_n \langle A_nh,k\rangle e_n.
\]
Since each $\langle A_nh,\cdot\rangle\in \mathfrak{g}_{CM}^*$ has a measurable linear
extension to $G$ such that $\|\langle
A_nh,\cdot\rangle\|_{L^2(\mu)} = \|A_nh\|_{\mathfrak{g}_{CM}}$ 
(see, for example, Theorem 2.10.11 of \cite{Bogachev1998}), we may extend
$\mathrm{ad}_h$ to a measurable linear transformation from $G=\mathfrak{g}$
to $G_{CM}=\mathfrak{g}_{CM}$ (still denoted by $\mathrm{ad}_h$) given by
\[ 
\mathrm{ad}_h g := \sum_n \langle A_nh,g\rangle e_n. \]
Note that here we are using the fact that
\begin{align*} \sum_n \|\langle A_nh,\cdot\rangle\|_{L^2(\mu)}^2
	= \sum_n \| A_nh\|_{H}^2
	&= \sum_{n,m} \langle A_nh,e_m\rangle^2 \\
	&= \sum_{n,m} \langle \mathrm{ad}_{h}e_m,e_n\rangle^2 
	\le \|h\|^2 \|[\cdot,\cdot]\|_{HS}^2
\end{align*}
which implies that
\[ \sum_n \langle A_nh,g\rangle ^2 < \infty \quad g\text{-a.s.} \]
Similarly, we may define
\[ \mathrm{ad}_g h := -\mathrm{ad}_h g
	= - \sum_n \langle A_nh,g\rangle e_n. \]

In a similar way, note that we may write, for $h,k\in G_{CM}$ and $m< r$,
\begin{align*} 
\mathrm{ad}_h^m k 
	&= \sum_{\ell_1}\cdots \sum_{\ell_m} \left(\prod_{b=1}^{m-1} 
		\langle A_{\ell_{b}}h, e_{\ell_{b+1}}\rangle  \right)\langle
		A_{\ell_m}h,k\rangle e_{\ell_1} \\
	&= (-1)^m \sum_{\ell_1}\cdots \sum_{\ell_m} \left(\prod_{b=1}^{m-1} 
		\langle A_{\ell_{b}}e_{\ell_{b+1}},h \rangle  \right)\langle
		A_{\ell_m}k,h\rangle e_{\ell_1}.
\end{align*}
and thus for $h,k,x\in G_{CM}$ and $n+m< r$
\begin{multline*} 
\mathrm{ad}_k^n \mathrm{ad}_h^m x 
	= (-1)^n \sum_{j_1}\cdots \sum_{j_n}
		\sum_{\ell_1}\cdots\sum_{\ell_m} 
		\left(\prod_{a=1}^{n-1} 
			\langle A_{j_{a}}e_{j_{a+1}},k \rangle  \right)
		\langle A_{j_n}e_{\ell_1},k\rangle \\
		\left(\prod_{b=1}^{m-1} 
			\langle A_{\ell_{b}}h, e_{\ell_{b+1}} \rangle  \right)
		\langle A_{\ell_m}h,x\rangle e_{j_1}.
\end{multline*}

More generally for $|n|+|m|<r$
\begin{align*} 
&\mathrm{ad}_k^{n_1} \mathrm{ad}_h^{m_1}\cdots 
		\mathrm{ad}_k^{n_s} \mathrm{ad}_h^{m_s} k \\
	&= (-1)^{|n|} \sum_{j^1_1}\cdots \sum_{\ell^s_{m_s}} 
		\prod_{c=1}^s \Bigg\{\left(\prod_{a_c=1}^{n_c-1} 
			\langle A_{j^c_{a_c}}e_{j^c_{a_c+1}},k \rangle  \right)
		\langle A_{j^c_{n_c}}e_{\ell^c_1},k\rangle \\
	&\qquad\qquad\times
		\left(\prod_{b_c=1}^{m_c-1} 
			\langle A_{\ell^c_{b_c}}h, e_{\ell^c_{b_c+1}} \rangle  \right)
		 \Bigg\}
		 \left\{\prod_{c=1}^{s-1}\langle A_{\ell^c_{m_c}}h,e_{j^{c+1}_1}\rangle\right\}
		\langle A_{\ell^s_{m_s}}h,k\rangle e_{j^1_1}\end{align*}
Thus, for $h\in G_{CM}$ and $|n|+|m|<r$, we may define a measurable action on $G$ given by
\begin{align*} 
G\ni g\mapsto\,&\mathrm{ad}_g^{n_1} \mathrm{ad}_h^{m_1}\cdots 
		\mathrm{ad}_g^{n_s} \mathrm{ad}_h^{m_s} g \\
	&:= (-1)^{|n|} \sum_{j^1_1}\cdots \sum_{\ell^s_{m_s}} 
		\prod_{c=1}^s \bigg\{\left(\prod_{a_c=1}^{n_c-1} 
			\langle A_{j^c_{a_c}}e_{j^c_{a_c+1}},g \rangle  \right)
		\langle A_{j^c_{n_c}}e_{\ell^c_1},g\rangle \\
	&\times
		\left(\prod_{b_c=1}^{m_c-1} 
			\langle A_{\ell^c_{b_c}}h, e_{\ell^c_{b_c+1}} \rangle  \right)
		 \Bigg\}
		 \left\{\prod_{c=1}^{s-1}\langle A_{\ell^c_{m_c}}h,e_{j^{c+1}_1}\rangle\right\}
		\langle A_{\ell^s_{m_s}}h,k\rangle e_{j^1_1}\in G_{CM}.
\end{align*}
(For $|n|+|m|\ge r$, we define this mapping to be 0, which is certainly
measurable.)
Again, we are using that 
\begin{multline*} 
\sum_{j^1_1}\Bigg\| \sum_{j^1_2}\cdots \sum_{\ell^s_{m_s}} 
		\prod_{c=1}^s \Bigg\{\left(\prod_{a_c=1}^{n_c-1} 
			\langle A_{j^c_{a_c}}e_{j^c_{a_c+1}},\cdot \rangle  \right)
		\langle A_{j^c_{n_c}}e_{\ell^c_1},\cdot\rangle \\
	\times
		\left(\prod_{b_c=1}^{m_c-1} 
			\langle A_{\ell^c_{b_c}}h, e_{\ell^c_{b_c+1}} \rangle  \right)
		 \Bigg\}
		 \left\{\prod_{c=1}^{s-1}\langle A_{\ell^c_{m_c}}h,e_{j^{c+1}_1}\rangle\right\}
		\langle A_{\ell^s_{m_s}}h,k\rangle e_{j^1_1} \Bigg\|_{L^2(\mu)}^2 <\infty.
\end{multline*}
This holds by straightforward but tedious computations --- in fact, iterative applications of
Cauchy-Schwarz combined with the fact that $\sum_{n,m} \|A_ne_m\|^2 =
\|[\cdot,\cdot]\|_{HS}^2<\infty$. Thus,  
\begin{multline*} 
\sum_{j^1_1}\Bigg(\sum_{j^1_2}\cdots \sum_{\ell^s_{m_s}} 
		\prod_{c=1}^s \Bigg\{\left(\prod_{a_c=1}^{n_c-1} 
			\langle A_{j^c_{a_c}}e_{j^c_{a_c+1}},g \rangle  \right)
		\langle A_{j^c_{n_c}}e_{\ell^c_1},g\rangle \\
	\times
		\left(\prod_{b_c=1}^{m_c-1} 
			\langle A_{\ell^c_{b_c}}h, e_{\ell^c_{b_c+1}} \rangle  \right)
		 \Bigg\}
		 \left\{\prod_{c=1}^{s-1}\langle A_{\ell^c_{m_c}}h,e_{j^{c+1}_1}\rangle\right\}
		\langle A_{\ell^s_{m_s}}h,k\rangle e_{j^1_1}\Bigg)^2 <\infty ,
\end{multline*}
$g$-a.s.~and $\mathrm{ad}_g^{n_1} \mathrm{ad}_h^{m_1}\cdots \mathrm{ad}_g^{n_s}
\mathrm{ad}_h^{m_s} g $ as given above is defined a.s.
Thus, we have the following result.

\begin{prop}
\label{p.gp}
For $h\in G_{CM}$, the mapping $G\ni g\mapsto g\cdot h\in G$ defined
analogously to (\ref{e.mult}) is a measurable right group action by
$G_{CM}$ on $G$, and similarly for the left action $g\mapsto h\cdot g$.
\end{prop}

\subsection{Examples of abstract nilpotent Lie algebras}

\begin{example}[Free nilpotent Lie algebras]
\label{ex.free}
Starting with an abstract Wiener space $(W,H,\mu)$, one may construct in the
standard way the abstract free
nilpotent Lie algebra of step $r$ with generators $\{h_i\}_{i=1}^\infty$ an
orthonormal basis of $H$.  See for example Section 0 of \cite{Gromov1996}.
\end{example}

\begin{example}[Heisenberg-like algebras]
\label{ex.heis}
Let $(W_i,H_i,\mu_i)$ for $i=1,2$ be abstract Wiener spaces.  Then for any
$\omega:H_1\times H_1\rightarrow H_2$ a Hilbert-Schmidt map, we may define a
Lie bracket on $\mathfrak{g}_{CM}=H_1\times H_2$ by
\[ [(h_1,h_2),(h_1',h_2')] := (0,\omega(h_1,h_1')), \]
and $\mathfrak{g}=W_1\times W_2$ may be thought of as an abstract
Heisenberg-like algebra as in \cite{DriverGordina2008}.
\end{example}

These abstract Heisenberg-like algebras are central extensions of
one abstract Wiener space by
another abstract Wiener space.  The next example generalizes this
construction.

\begin{example}[Extensions of Lie algebras]
\label{ex.ext}
Let $\mathfrak{v}$ and $\mathfrak{h}$ be Lie
algebras, and let $\mathrm{Der}(\mathfrak{v})$ denote the set of derivations
on $\mathfrak{v}$; that is, $\mathrm{Der}(\mathfrak{v})$ consists of all
linear maps $\rho:\mathfrak{v}\rightarrow\mathfrak{v}$ satisfying Leibniz's
rule:
\[ \rho([X,Y]_\mathfrak{v}) = [\rho(X),Y]_\mathfrak{v} +
	[X,\rho(Y)]_\mathfrak{v}. \]
Now suppose there are a linear mapping 
$\alpha: \mathfrak{h} \rightarrow \mathrm{Der}(\mathfrak{v})$ 
and a skew-symmetric bilinear mapping
$\omega:\mathfrak{h}\times \mathfrak{h}\rightarrow\mathfrak{v}$, 
satisfying, for all $X,Y,Z\in \mathfrak{h}$,
\begin{equation}
\label{b1} 
\tag{B1}
[\alpha_X,\alpha_Y] - \alpha_{[X,Y]_\mathfrak{h}} 
	= \mathrm{ad}_{\omega(X,Y)} 
\end{equation}
and
\begin{equation}
\label{b2} 
\tag{B2}
\sum_{\text{cyclic}}\left(\alpha_X \omega(Y,Z) -
	\omega([X,Y]_\mathfrak{h},Z)\right)= 0.  
\end{equation}
Then, one may verify that, for $X_1+V_1,X_2+V_2\in \mathfrak{h}\oplus
\mathfrak{v}$,
\[ [X_1+V_1, X_2+V_2]_\mathfrak{g} :=  [X_1,X_2]_\mathfrak{h} + \omega(X_1,X_2) 
	+ \alpha_{X_1}V_2 -\alpha_{X_2}V_1 + [V_1,V_2]_\mathfrak{v} \]
defines a Lie bracket on $\mathfrak{g}:=\mathfrak{h}\oplus\mathfrak{v}$, 
and we say $\mathfrak{g}$ is an extension of $\mathfrak{h}$ over 
$\mathfrak{v}$. 
That is, $\mathfrak{g}$ is the Lie algebra with ideal
$\mathfrak{v}$ and quotient algebra
$\mathfrak{g}/\mathfrak{v}=\mathfrak{h}$.  The associated exact sequence is
\[ 0\rightarrow\mathfrak{v}\overset{\iota_1}{\longrightarrow}\mathfrak{g}
	\overset{\pi_2}{\longrightarrow}\mathfrak{h}\rightarrow0,\] 
where $\iota_1$ is
inclusion and $\pi_2$ is projection.  In fact,  these are the 
only extensions of
$\mathfrak{h}$ over $\mathfrak{v}$
(see, for example, \cite{AMR00}).

Now suppose that $(W,H,\mu)$ is a real abstract Wiener space, and
$(\mathfrak{v},\mathfrak{v}_{CM},\mu^0)$ is an abstract nilpotent Lie algebra.
Motivated by the previous discussion, we may
consider $H$ as an abelian Lie algebra and construct extensions of
$H$ over $\mathfrak{v}_{CM}$.
In this case, we need  a linear mapping
$\alpha:H\rightarrow \mathrm{Der}(\mathfrak{v}_{CM})$
and a skew-symmetric bilinear mapping
$\omega:H\times H\rightarrow\mathfrak{v}_{CM}$,
such that $\omega$ and
$\alpha:H\times\mathfrak{v}_{CM}\rightarrow\mathfrak{v}_{CM}$ 
are both Hilbert-Schmidt and together $\omega$ and $\alpha$
satisfy (\ref{b1}) and (\ref{b2}), which in
this setting become
\begin{equation*}
[\alpha_X,\alpha_Y] = \mathrm{ad}_{\omega(X,Y)}
\end{equation*}
and
\begin{equation*}
\alpha_X \omega(Y,Z) + \alpha_Y \omega(Z,X) + \alpha_Z \omega(X,Y) = 0, 
\end{equation*}
for all $X,Y,Z\in H$.  Then we may define a Lie algebra structure
on $\mathfrak{g}_{CM}= H\oplus \mathfrak{v}_{CM}$ via the Lie bracket
\begin{equation*} 
[(X_1,V_1), (X_2,V_2)]_{\mathfrak{g}_{CM}} := (0, \omega(X_1,X_2)
	+ \alpha_{X_1}V_2 - \alpha_{X_2}V_1 ). 
\end{equation*}
Again, these are in fact the only extensions of $H$
over $\mathfrak{v}_{CM}$.  
\end{example}

Combining Examples \ref{ex.heis} and \ref{ex.ext} shows that these
constructions are iterative, in that one may construct new abstract nilpotent
Lie algebras as Lie algebra extensions of another abstract nilpotent Lie
algebra.  The next example builds on the previous one to give one precise way
to construct some Lie algebra extensions.

\begin{example}
\label{ex.beta}
Let $\beta:H\rightarrow \mathfrak{v}_{CM}$ be any Hilbert-Schmidt map, and for
$X,Y\in H$ and $V\in \mathfrak{v}_{CM}$ define
\[ \omega(X,Y)  := [\beta(X),\beta(Y)]_{\mathfrak{v}_{CM}} \] 
\[ \alpha_X V := \mathrm{ad}_{\beta(X)} V = [\beta(X),V]_{\mathfrak{v}_{CM}}. \] 
Then the Lie bracket on $\mathfrak{g}_{CM}$ is given by
\begin{multline*}
[(X,V),(Y,U)]  \\
	= (0,[\beta(X),\beta(Y)]_{\mathfrak{v}_{CM}} 
	+ [\beta(X),U]_{\mathfrak{v}_{CM}} - [\beta(Y),V]_{\mathfrak{v}_{CM}} 
	+ [V,U]_{\mathfrak{v}_{CM}}). 
\end{multline*}
Note that, if $\mathfrak{v}_{CM}$ is nilpotent of step $r$, then
$\mathfrak{g}_{CM}$ will automatically be nilpotent of step $r$. 
\end{example}

The examples above demonstrate that the space of abstract nilpotent Lie
algebras is quite rich, and there are many natural examples with a straightforward construction.  This
significantly improves the results of \cite{Melcher2009}, which studied 
heat kernel measures on nilpotent extensions of abstract
Wiener spaces over {\it finite-dimensional} nilpotent Lie algebras.
For example, the restriction $\mathrm{dim}(\mathfrak{v}_{CM})<\infty$
trivializes Example \ref{ex.beta}.  We elaborate in the following remark.

\begin{remark}
Returning to Example \ref{ex.beta}, since $\beta$ is linear and continuous, we have the decomposition
\[ H = \mathrm{Nul}(\beta) \oplus \mathrm{Nul}(\beta)^\perp, \] 
where $\mathrm{dim}(\mathrm{Nul}(\beta)^\perp) \le
\mathrm{dim}(\mathfrak{v}_{CM})$.  Thus, for $X,Y\in H$ we may write
$X=X_1+X_2, Y=Y_1+Y_2\in \mathrm{Nul}(\beta) \oplus \mathrm{Nul}(\beta)^\perp$
and
\begin{align*} 
[(X_1+X_2,0),(Y_1+Y_2,0)] 
	&= [ \beta(X_1+X_2),\beta(Y_1+Y_2)] 
	= [\beta(X_2),\beta(Y_2)], 
\end{align*}
and $\omega$ is a map on $\mathrm{Nul}(\beta)^\perp\times
\mathrm{Nul}(\beta)^\perp$.  
Thus, $[\mathrm{Nul}(\beta),\mathrm{Nul}(\beta)]=\{0\}$ and similarly
$[\mathrm{Nul}(\beta),\mathfrak{v}]=\{0\}$.  So
\[ \mathfrak{g}_{CM}= H\oplus\mathfrak{v}_{CM} = \mathrm{Nul}(\beta) \oplus
	\mathrm{Nul}(\beta)^\perp\oplus\mathfrak{v}_{CM}, \]
and, in particular, when $\mathrm{dim}(\mathfrak{v}_{CM})<\infty$, $\mathfrak{g}_{CM}$ 
is just an extension of the finite-dimensional vector space
$\mathrm{Nul}(\beta)^\perp$ by the finite-dimensional Lie algebra
$\mathfrak{v}_{CM}$, and the construction is not truly infinite-dimensional.
\end{remark}

\subsection{Properties of $\mathfrak{g}_{CM}$}

This section collects some results for
topological and geometric properties of $\mathfrak{g}_{CM}$ that we'll require for
the sequel.  

\begin{prop}
\label{p.normest}
For all $m\ge2$, $[[[\cdot,\cdot],\ldots],\cdot]:\mathfrak{g}_{CM}^{\otimes m}
\rightarrow \mathfrak{g}_{CM}$ is Hilbert-Schmidt.
\end{prop}

\begin{proof}
For $m=2$, this follows from the definition of $G$.  Now assume the statement holds for all $m=\ell$, and 
consider $m=\ell+1$.  Writing $[[h_{i_1},h_{i_2}],\cdots,h_{i_\ell+1}]\in 
\mathfrak{g}_{CM}$ in terms of the orthonormal basis
$\{e_j\}_{j=1}^\infty$ and using 
H\"older's inequality gives
\begin{align*}
\|[[[\cdot,&\cdot],\ldots],\cdot]\|_2^2
	= \|[[[\cdot,\cdot],\ldots],\cdot]\|_{(\mathfrak{g}_{CM}^*)^{\otimes \ell+1}
		\otimes\mathfrak{g}_{CM}} \\
	&= \sum_{i_1,\ldots,i_{\ell+1}=1}^\infty 
		\|[[[h_{i_1},h_{i_2}],\cdots,h_{i_\ell}],h_{i_{\ell+1}}]\|
		^2 \\
	&= \sum_{i_1,\ldots,i_{\ell+1}=1}^\infty \left\|\sum_{j=1}^\infty 
		[e_j,h_{i_{\ell+1}}]
		\langle e_j, [[h_{i_1},h_{i_2}],\cdots,h_{i_\ell}]\rangle
		\right\|^2 \\
	&\le \sum_{i_1,\dots,i_{\ell+1}=1}^\infty \left(\sum_{j=1}^\infty 
			\| [e_j,h_{i_{\ell+1}}] \|^2\right)
		\left(\sum_{j=1}^\infty 
		|\langle e_j,[[h_{i_1},h_{i_2}],\cdots,h_{i_\ell}]\rangle|^2\right) \\
	&= \left(\sum_{i_{\ell+1}=1}^\infty \sum_{j=1}^\infty 
			\| [e_j,h_{i_{\ell+1}}] \|^2\right)
		\left(\sum_{i_1,\dots,i_\ell=1}^\infty 
		\|[[h_{i_1},h_{i_2}],\cdots,h_{i_\ell}]\|^2\right)
		\\
	&\le \|[\cdot,\cdot]\|_{(\mathfrak{g}_{CM}^*)^{\otimes 2}
		\otimes\mathfrak{g}_{CM}}^2 \|[[[\cdot,\cdot],\ldots],\cdot]\|_{(\mathfrak{g}_{CM}^*)^{\otimes \ell}
		\otimes\mathfrak{g}_{CM}}^2
\end{align*}
and the last line is finite by the induction hypothesis.
\end{proof}

Next we will recall that the flat and geometric topologies on
$\mathfrak{g}_{CM}$ are equivalent.  First we set the following notation.

\begin{notation}
For $g\in G_{CM},$ let $L_g:G_{CM}\rightarrow G_{CM}$ and $R_g:G_{CM}\rightarrow G_{CM}$
denote left and right multiplication by $g$, respectively.  As $G_{CM}$
is a vector space, to each $g\in G_{CM}$ we can associate the tangent
	space $T_g G_{CM}$ to $G_{CM}$ at $g$, which is naturally isomorphic to $G_{CM}$.
For $f:G_{CM}\rightarrow\mathbb{R}$ a Frech\'{e}t smooth function and  
	$v,x\in G_{CM}$ and $h\in\mathfrak{g}_{CM}$, let
\[ f'(x)h := \partial_h f(x) = \frac{d}{dt}\bigg|_{0}f(x+th), \]
	and let $v_x \in T_x G_{CM}$ denote the tangent vector
satisfying $v_xf=f'(x)v$.  If $\sigma(t)$ is any smooth curve in
$G_{CM}$ such that $\sigma(0) = x$ and $\dot{\sigma}(0)=v$ (for example,
$\sigma(t) = x+tv$), then
\[ L_{g*} v_x = \frac{d}{dt}\bigg|_0 g\cdot \sigma(t). \]
\end{notation}

\begin{notation}
\label{n.length}
Let $C^1([0,1],G_{CM})$ denote the collection of $C^1$-paths
$g:[0,1]\rightarrow G_{CM}$.  The length of $g$ is defined as 
\[ \ell_{CM}(g) := \int_0^T \|L_{g^{-1}(s)*}g'(s)\|_{\mathfrak{g}_{CM}}\,ds.
\]
The Riemannian distance between $x,y\in G_{CM}$ is then defined as
\begin{align*} d(x,y) := \inf\{\ell_{CM}(g): g\in C_h^1([0,1],G_{CM}) \text{ such that }
	g(0)=x \text{ and } g(1)=y \}. 
\end{align*}
\end{notation}

The following proposition was proved as Corollary 4.13 in
\cite{Melcher2009}.  This was proved under the conditions that
$\mathfrak{g}_{CM}$ was a ``semi-infinite Lie algebra'', that is,
under the assumption that $\mathfrak{g}_{CM}$ was a nilpotent Lie algebra
extension (as in Example \ref{ex.ext}) 
of a finite-dimensional nilpotent Lie algebra $\mathfrak{v}$ by an
abstract Wiener space.  However, a cursory inspection of the proofs
there will show that they only depended on the fact that the Lie
bracket was Hilbert-Schmidt, and not on the ``stratified'' structure
of $\mathfrak{g}_{CM}=H\times\mathfrak{v}$ or the fact that the image of the Lie bracket
was a finite-dimensional subspace.

\begin{prop}
\label{p.length}
The topology on $G_{CM}$ induced by $d$ is equivalent to the Hilbert
topology induced by $\|\cdot\|_{\mathfrak{g}_{CM}}$.
\end{prop}

We may find a uniform lower bound on the Ricci curvature of all
finite-dimensional subgroups of $G_{CM}$.  Finite-dimensional subgroups
of $G_{CM}$ may be obtained by taking the Lie algebra
generated by any finite-dimensional subspace of $\mathfrak{g}_{CM}$
(which is again necessarily finite-dimensional by the nilpotence
of the bracket) and endowing it with the standard group operation via the
Baker-Campbell-Hausdorff-Dynkin formula.   An analogue of the following proposition was
proved as Proposition 3.23 and Corollary 3.24 in \cite{Melcher2009}.
It follows directly from the form of the Ricci curvature on nilpotent
groups (endowed with a left invariant metric) and the assumption
that the Lie bracket is Hilbert-Schmidt.  This proof is essentially the same as in
\cite{Melcher2009}, but it is quite brief and so is included for completeness.

\begin{prop}
\label{p.Ric}
Let 
\[ k := -\frac{1}{2} \sup \left\{ 
		\|[\cdot,X]\|^2_{\mathfrak{g}_{CM}^*\otimes\mathfrak{g}_{CM}}   
		: \,\|X\|_{\mathfrak{g}_{CM}}=1 \right\}. \]
Then $k>-\infty$ and $k$ is the largest constant such that
\[ \langle \mathrm{Ric}^\pi X,X\rangle_{\mathfrak{g}_\pi} \ge 
	k \|X\|^2_{\mathfrak{g}_\pi}, 
	\quad \text{ for all } X \in\mathfrak{g}_\pi, \]
holds uniformly for all $\mathfrak{g}_\pi$ finite-dimensional Lie subalgebras
of $\mathfrak{g}_{CM}$.
\end{prop}

\begin{proof}
For $\mathfrak{g}$ any nilpotent Lie algebra with orthonormal basis
$\Gamma$, 
\begin{align*}
\langle \mathrm{Ric}\, X,X\rangle 
	&= \frac{1}{4}\sum_{Y\in\Gamma} \|\mathrm{ad}^*_Y X\|^2 
		- \frac{1}{2}\sum_{Y\in\Gamma} \|\mathrm{ad}_Y X\|^2 
	\ge - \frac{1}{2}\sum_{Y\in\Gamma} \|[Y,X]\|^2
\end{align*}
for all $X\in\mathfrak{g}$.  Thus, for $\mathfrak{g}_\pi$ any
finite-dimensional Lie algebra
\[ \langle \mathrm{Ric}^\pi X,X\rangle_{\mathfrak{g}_\pi} \ge 
	k_\pi \|X\|^2_{\mathfrak{g}_\pi}, 
	\quad \text{ for all } X \in\mathfrak{g}_\pi, \]
where
\begin{equation}
\label{e.pah}
 k_\pi :=  - \frac{1}{2}\sup \left\{ 
		\|[\cdot,X]\|^2_{\mathfrak{g}_\pi^*\otimes\mathfrak{g}_\pi}
 		:\, \|X\|_{\mathfrak{g}_\pi}=1 \right\} 
	\ge - \frac{1}{2}\|[\cdot,\cdot]\|_2^2 > -\infty.
\end{equation}
Taking the infimum of $k_\pi$ over all $\mathfrak{g}_\pi$ completes the proof.
\end{proof}

\section{Brownian motion on $G$}
\label{s.BM}
Suppose that $B_t$ is a smooth curve in $\mathfrak{g}_{CM}$ with
$B_0=0$, and consider the differential equation
\[ \dot{g}_t = g_t \dot{B}_t 
	:= L_{g_t*}\dot{B}_t, \quad \text{ with } g_0=\mathbf{e}. \]  
The solution $g_t$ may be written as follows (see
\cite{Strichartz87}):  For $t>0$, let $\Delta_n(t)$ denote the simplex in
$\mathbb{R}^n$ given by
\[ \{s=(s_1,\cdots,s_n)\in\mathbb{R}^n: 0<s_1<s_2<\cdots<s_n<t\}. \]
Let $\mathcal{S}_n$ denote the permutation group on $(1,\cdots,n)$, and, for each
$\sigma\in \mathcal{S}_n$, let $e(\sigma)$ denote the number of ``errors'' in the
ordering $(\sigma(1),\sigma(2),\cdots,\sigma(n))$, that is, $e(\sigma)=\#
\{j<n: \sigma(j)>\sigma(j+1)\}$.  Then
\begin{multline}
\label{e.ode}
g_t = \sum_{n=1}^r \sum_{\sigma\in \mathcal{S}_n} 
	\left( (-1)^{e(\sigma)}\bigg/ n^2 
	\begin{bmatrix} n-1 \\ e(\sigma) \end{bmatrix}\right) \times \\
	\int_{\Delta_n(t)} [ \cdots[\dot{B}_{s_{\sigma(1)}},
	\dot{B}_{s_{\sigma(2)}}],\ldots,] \dot{B}_{s_{\sigma(n)}}]\, ds,
\end{multline}
where the $n=1$ term is understood to be $\int_0^t dB_s = B_t$.
Using this as our motivation, we first explore stochastic integral
analogues of equation (\ref{e.ode}) where the smooth curve $B$ is replaced by
Brownian motion on $\mathfrak{g}$.

\subsection{Brownian motion and finite-dimensional approximations}

We now return to the setting of an abstract Wiener space
$(\mathfrak{g},\mathfrak{g}_{CM},\mu)$ endowed with a nilpotent
Hilbert-Schmidt Lie bracket on $\mathfrak{g}_{CM}$.
Again, let $B_t$ denote Brownian motion on $\mathfrak{g}$.
By equation (\ref{e.ode}), the solution to the Stratonovich
stochastic differential equation
\[ \delta g_t = L_{g_t*} \delta B_t, \quad \text{ with } g_0=\mathbf{e}, \]
should be given by
\begin{equation}
\label{e.gt}
g_t = \sum_{n=1}^{r} \sum_{\sigma\in\mathcal{S}_n} c^\sigma_n \int_{\Delta_n(t)} 
	[ [\cdots[\delta B_{s_{\sigma(1)}},\delta B_{s_{\sigma(2)}}],\cdots], 
	\delta B_{s_{\sigma(n)}}],
\end{equation}
for coefficients $c_n^\sigma$ determined by equation (\ref{e.ode}).

To understand the integrals in (\ref{e.gt}), consider the following heuristic
computation.
Let $\{M_n(t)\}_{t\ge0}$ denote the process in $\mathfrak{g}^{\otimes n}$ 
defined by
\[ 
M_n(t) := \int_{\Delta_n(t)} \delta B_{s_1}\otimes\cdots\otimes \delta
	B_{s_n}.
\]
By repeatedly applying the definition of the Stratonovich integral, 
the iterated Stratonovich integral $M_n(t)$ 
may be realized as a linear combination of iterated It\^o integrals:
\[ M_n(t) = \sum_{m=\lceil n/2\rceil}^n \frac{1}{2^{n-m}}
		\sum_{\alpha\in\mathcal{J}_n^m} I^n_t(\alpha), \]
where
\[ \mathcal{J}_n^m := \left\{(\alpha_1,\ldots,\alpha_m)\in\{1,2\}^m :
	\sum_{i=1}^m \alpha_i = n \right\}, \]
and, for $\alpha\in\mathcal{J}_n^m$, $I_t^n(\alpha)$ is the iterated 
It\^o integral
\[ I_t^n(\alpha) = \int_{\Delta_m(t)} dX^1_{s_1}\otimes\cdots\otimes 
	dX^m_{s_m} \]
with
\[ dX_s^i = \left\{ \begin{array}{cl} dB_s & \text{if } \alpha_i=1 \\
	\sum_{j=1}^\infty h_j \otimes h_j \, ds & \text{if } \alpha_i=2
	\end{array} \right. ; \]
compare with Proposition 1 of \cite{BenArous89}.  This change from multiple
Stratonovich integrals to multiple It\^o integrals may also be recognized as a
specific case of the Hu-Meyer formulas \cite{HuMeyer88-2,HuMeyer88-1},
but we will compute more explicitly to verify that our integrals are well-defined.

Define $F_1:\mathfrak{g}_{CM}\to\mathfrak{g}_{CM}$ by $F_1(k)=k$, and for $n\in\{2,\cdots,r\}$ define $F_n:\mathfrak{g}_{CM}^{\otimes n}
\rightarrow\mathfrak{g}_{CM}$ by
\begin{equation}
\label{e.1Fn}
 F_n(k_1\otimes\cdots\otimes k_n) 
	:= [ [ [\cdots[k_1,k_2],k_3],\cdots],k_n]. 
\end{equation}
For each fixed $n$ and $\sigma\in\mathcal{S}_n$, define $F_n^\sigma:\mathfrak{g}_{CM}^{\otimes n}
\rightarrow\mathfrak{g}_{CM}$ by
\begin{equation}
\label{e.Fn}
\begin{split}
F_n^\sigma(k_1\otimes\cdots\otimes k_n) 
	&:= F_n(k_{\sigma(1)}\otimes\cdots\otimes k_{\sigma(n)}) \\
	&= [ [\cdots[k_{\sigma(1)},k_{\sigma(2)}],\cdots], 
		k_{\sigma(n)}].
\end{split}
\end{equation}
Then we may write
\[
g_t = \sum_{n=1}^{r} \sum_{\sigma\in\mathcal{S}_n} 
	c^\sigma_n F^\sigma_n (M_n(t))
	= \sum_{n=1}^{r} \sum_{\sigma\in\mathcal{S}_n}
	\sum_{m=\lceil n/2\rceil}^n \frac{c^\sigma_n }{2^{n-m}}
		\sum_{\alpha\in\mathcal{J}_n^m} F^\sigma_n (I^n_t(\alpha)),
\]
presuming we can make sense of the integrals $F_n^\sigma(I_t^n(\alpha))$.

For each $\alpha$, let $p_\alpha=\#\{i:\alpha_i=1\}$ and $q_\alpha=\#\{i:
\alpha_i=2\}$ (so that $p_\alpha+q_\alpha=m$ when $\alpha\in\mathcal{J}_n^m$),
and let 
\[ \mathcal{J}_n := \bigcup_{m=\lceil n/2\rceil}^n \mathcal{J}_n^m. \]
Then, for each $\sigma\in\mathcal{S}_n$ and $\alpha\in\mathcal{J}_n$, 
\[ F_n^\sigma(I_t^n(\alpha))
	= \int_{\Delta_{p_\alpha}(t)} f_\alpha(s,t) \hat{F}_n^{\sigma,\alpha}
		(dB_{s_1}\otimes\cdots\otimes dB_{s_{p_\alpha}}),
\]
where $\hat{F}_n^{\sigma,\alpha}$ and $f_\alpha$ are as follows.

The map $\hat{F}_n^{\sigma,\alpha}:\mathfrak{g}^{\otimes p_\alpha}
\rightarrow\mathfrak{g}$ is defined by
\begin{multline}
\label{e.Fhat}
\hat{F}_n^{\sigma,\alpha}(k_1\otimes\cdots\otimes k_{p_\alpha}) \\
	:= \sum_{j_1,\ldots,j_{q_\alpha}=1}^\infty 
		F_n^{\sigma'}(k_1\otimes\cdots\otimes k_{p_\alpha}
		\otimes h_{j_1}\otimes h_{j_1}
		\otimes\cdots\otimes h_{j_{q_\alpha}}\otimes h_{j_{q_\alpha}}), 
\end{multline}
for $\{h_j\}_{j=1}^\infty$ an orthonormal basis of $\mathfrak{g}_{CM}$ and
$\sigma'=\sigma'(\alpha)\in\mathcal{S}_n$ given by 
$\sigma'=\sigma\circ\tau^{-1}$, for any $\tau\in\mathcal{S}_n$ such that
\begin{multline*} 
\tau(dX^1_{s_1}\otimes\cdots\otimes dX^m_{s_m}) \\
	= \sum_{j_1,\cdots,j_{q_\alpha}=1}^\infty dB_{s_1}\otimes\cdots
	\otimes dB_{s_{p_\alpha}}\otimes h_{j_1}\otimes h_{j_1}\otimes\cdots
	\otimes h_{j_{q_\alpha}}\otimes h_{j_{q_\alpha}} ds_1\cdots
	ds_{q_\alpha}.
\end{multline*}

To define $f_\alpha$, first consider the polynomial of order $q_\alpha$,
in the variables $s_i$ with $i$ such that $\alpha_i=1$ and in the variable $t$, given by evaluating the integral
\begin{equation}
\label{e.fprime} 
f_\alpha'( (s_i:\alpha_i=1),t)
	= \int_{\Delta'_{q_\alpha}(t)} \prod_{i: \alpha_i=2} ds_i,
\end{equation}
where $\Delta'_{q_\alpha}(t)=\{s_{i-1}<s_i<s_{i+1}: \alpha_i=2\}$ with 
$s_0=0$ and $s_{m+1}=t$.  
Then $f_\alpha$ is $f_\alpha'$ with the variables reindexed by the
bijection $\{i: \alpha_i=1\}\rightarrow\{1,\ldots,p_\alpha\}$ that
maintains the natural ordering of these sets.
(For example, for $\alpha=(1,2,1,2)\in\mathcal{J}_6^4$,
\[ f_\alpha'(s_1,s_3,t)=\int_{\{s_1<s_2<s_3,s_3<s_4<t\}} ds_2\,ds_4
	= (t-s_3)(s_3-s_1) , \]
so that $f_\alpha(s_1,s_2,t)=(t-s_2)(s_2-s_1)$.) 

This explicit realization of $f_\alpha$ is not critical to the sequel.
It is really only necessary to know that 
$f_\alpha$ is a polynomial of order $q_\alpha$ in
$s=(s_1,\ldots,s_{p_\alpha})$ and $t$, and thus may be written as
\[ f_\alpha(s,t) = \sum_{a=0}^{q_\alpha} b_\alpha^a t^a 
	\tilde{f}_{\alpha,a}(s), 
\]
for some coefficients $b_\alpha^a\in\mathbb{R}$ and polynomials 
$\tilde{f}_{\alpha,a}$ of degree $q_\alpha-a$ in $s$.  
Now, if $\hat{F}_n^{\sigma,\alpha}$ is
Hilbert-Schmidt on $\mathfrak{g}_{CM}^{\otimes p_\alpha}$, then
\[ \int_{\Delta_{p_\alpha}(t)} \left\|\tilde{f}_{\alpha,a}(s)
	\hat{F}_n^{\sigma,\alpha} \right\|_2^2\, ds	
	= \left\|\tilde{f}_{\alpha,a}\right\|^2_{L^2(\Delta_{p_\alpha}(t))}
		\left\|\hat{F}_n^{\sigma,\alpha} \right\|_2^2
	<\infty,
\]
and 
\begin{align}
\label{e.a}
F_n^\sigma(I_t^n(\alpha))
	&= \sum_{a=0}^{q_\alpha} b_\alpha^a t^a J_n(\tilde{f}_{\alpha,a}
		\hat{F}_n^{\sigma,\alpha})_t
\end{align}
may be understood in the sense of the limit integrals in 
Proposition \ref{p.int2}.  (In particular, if $\alpha_m=1$, then 
$f_\alpha=f_\alpha(s)$ does not depend on $t$, and Proposition \ref{p.int2} 
implies that $F_n^\sigma(I_t^n(\alpha))$ is a $\mathfrak{v}$-valued 
$L^2$-martingale.)

The above computations show that, if for all $n$, $\sigma\in\mathcal{S}_n$,
and $\alpha\in\mathcal{J}_n$,
$\hat{F}_n^{\sigma,\alpha}$ is Hilbert-Schmidt, then we may rewrite (\ref{e.gt}) as
\[ g_t = \sum_{n=1}^{r} \sum_{\sigma\in\mathcal{S}_n}
	\sum_{m=\lceil n/2\rceil}^n \frac{c^\sigma_n }{2^{n-m}}
		\sum_{\alpha\in\mathcal{J}_n^m}  
	\sum_{a=0}^{q_\alpha} b_\alpha^a t^a J_n(\tilde{f}_{\alpha,a}
		\hat{F}_n^{\sigma,\alpha})_t,
\]
where $J_n$ is as defined in Proposition \ref{p.int2}.  
The next proposition shows that  $\hat{F}_n^{\sigma,\alpha}$ is 
Hilbert-Schmidt as desired, and thus $g_t$ in (\ref{e.gt}) is well-defined.

\begin{prop}
\label{p.HS}
Let $n\in\{2,\ldots,r\}$, $\sigma\in\mathcal{S}_n$, and 
$\alpha\in\mathcal{J}_n$.
Then $\hat{F}_n^{\sigma,\alpha}:\mathfrak{g}_{CM}^{\otimes p_\alpha}
\rightarrow \mathfrak{g}_{CM}$ is Hilbert-Schmidt.
\end{prop}

\begin{proof}
For the whole of this proof, all sums will be taken over an orthonormal basis
of $\mathfrak{g}_{CM}$.

Now, for $F_n$ and $F_n^\sigma$ as defined in equations (\ref{e.1Fn}) and
(\ref{e.Fn}), we may write
\begin{align*} 
F_n^\sigma(k_1\otimes\cdots\otimes k_{p_\alpha}\otimes h_1\otimes
		h_1\otimes\cdots\otimes h_{q_\alpha}\otimes h_{q_\alpha}) 
	= F_n(A_{\sigma(1)}\otimes\ldots\otimes A_{\sigma(n)})
\end{align*}
where
\[ 
A_b := \left\{\begin{array}{ll}
		k_b & \text{if } b=1,\ldots,p_\alpha \\
		h_{\lceil (b-p_\alpha)/2\rceil} & \text{if } b=p_\alpha + 1,\ldots,n
	\end{array} \right. .
\] 

If $q_\alpha=0$, then $p_\alpha=n$, each $A_{\sigma(j)}=k_i$ for some
$i=1,\ldots,n$, and
\begin{align*}
\|\hat{F}_n^{\sigma,\alpha}\|_2^2
	&= \sum_{k_1,\ldots,k_n} 
		\|F_n^\sigma(k_1\otimes\cdots\otimes k_n)\|_{\mathfrak{g}_{CM}}^2 \\
	&= \sum_{k_1,\ldots,k_n} 
		\|F_n(k_1\otimes\cdots\otimes k_n)\|_{\mathfrak{g}_{CM}}^2
	= \|[ [[\ldots[\cdot,\cdot],\ldots],\cdot]\|_{(\mathfrak{g}_{CM}^*)^{\otimes
	n}\otimes\mathfrak{g}_{CM}}^2
\end{align*}
which is finite by Proposition \ref{p.normest}.

Now, if $q_\alpha=1$, let $N=N(\sigma)$ denote the second
$j$ such that $\sigma(j)\in\{n-1,n\}$; that is,
\begin{align*}
F_n^\sigma(k_1\otimes\cdots\otimes& k_{n-1}\otimes h_1\otimes
		h_1) \\
	&= F_n(A_{\sigma(1)}\otimes \cdots\otimes A_{\sigma(N-1)}\otimes h_1 \otimes
		A_{\sigma(N+1)}\otimes\ldots\otimes A_{\sigma(n)}) \\
	&= [[[[ [ \ldots[A_{\sigma(1)},A_{\sigma(2)}],\ldots],A_{\sigma(N-1)}]
		,h_1],A_{\sigma(N+1)}],\ldots],A_{\sigma(n)}]
\end{align*} 
where exactly one of $A_{\sigma(1)},\ldots,A_{\sigma(N-1)}$ is $h_1$ and,
for all $j>N$, $A_{\sigma(j)}=k_i$ for some 
\[ i\in I := I(\sigma):= \{ \sigma(j):j=N+1,\ldots,n\} \subseteq
 \{1,\ldots,p_\alpha\}= \{1,\ldots,n-2\}. \]
Thus, writing
\[ \mathcal{A}(h_1, k_i : i\in I^c) 
	:= [ [ \ldots[A_{\sigma(1)},A_{\sigma(2)}],\ldots],A_{\sigma(N-1)}], \]
we have that
\begin{multline*}
F_n^\sigma(k_1\otimes\cdots\otimes k_{n-1}\otimes h_1\otimes h_1) \\
	= \sum_{e_1}
		\,\langle \mathcal{A}(h_1,k_i:i\in I^c),e_1\rangle_{\mathfrak{g}_{CM}}
		[ [ \ldots[e_1,h_1],A_{\sigma(N+1)}],\ldots,A_{\sigma(n)}],
\end{multline*}
and so
\begin{align*}
\|\hat{F}_n^{\sigma,\alpha}\|_2^2
	&= \sum_{k_1,\ldots,k_{n-1}}\left\|
		\sum_{h_1} \,F_n^\sigma(k_1\otimes\cdots\otimes k_{n-1}\otimes h_1\otimes h_1)
		\right\|_{\mathfrak{g}_{CM}}^2 \\
	&\le \sum_{k_1,\ldots,k_{n-1}} \left( \sum_{h_1,e_1} 
		|\langle \mathcal{A}(h_1,k_i :i\in I^c),e_1\rangle
		_{\mathfrak{g}_{CM}}|^2 \right) \\
	&\qquad\qquad\qquad\quad\times
		\left( \sum_{h_1,e_1} 
		\|[ [[e_1,h_1],A_{\sigma(N+1)}],\ldots,A_{\sigma(n)}]\|
		_{\mathfrak{g}_{CM}}^2 \right) \\
	&= \left(\sum_{k_i:i\in I^c,h_1,e_1} |\langle \mathcal{A}(h_1,k_i :i\in I^c),e_1\rangle
		_{\mathfrak{g}_{CM}}|^2 \right) \\
	&\qquad\qquad\qquad\quad \times 
		\left( \sum_{k_i:i\in I,h_1,e_1} 
		\|[ [\ldots[[e_1,h_1],A_{\sigma(N+1)}],\ldots],A_{\sigma(n)}]\|
		_{\mathfrak{g}_{CM}}^2 \right) \\
	&\le \|[ [ \ldots[\cdot,\cdot],\ldots],\cdot]\|
		_{(\mathfrak{g}_{CM}^*)^{\otimes N-1}\otimes\mathfrak{g}_{CM}}^2
		\|[ [ \ldots[\cdot,\cdot],\ldots],\cdot]\|
		_{(\mathfrak{g}_{CM}^*)^{\otimes n-N+1}\otimes\mathfrak{g}_{CM}}^2.
\end{align*}
 
Now more generally when $q_\alpha\ge2$, we may similarly ``separate'' the
pairs of $h_j$'s as above.  More precisely, define
\[ \Phi(b) := \Phi_\alpha(b) := \left\{  \begin{array}{ll} 
	b & \text{if } b=1,\ldots, p_\alpha \\
	\lceil\frac{b-p_\alpha}{2}\rceil+p_\alpha & \text{if } b=p_\alpha+1,\ldots,n
\end{array}\right. .\]
Let $N_0=1$, and set
\[ \Omega_1^j := \{ \Phi(\sigma(\ell)) : \ell=N_0,\ldots,j-1\} \quad
\text{ and } \quad N_1 := \min\{ j>N_0:
	\Phi(\sigma(j))\in\Omega_1^j\},  \]
\[ \Omega_2^j := \{\Phi(\sigma(\ell)):\ell=N_1,\ldots,j-1\} \quad
\text{ and } \quad N_2 := \min\{ j>N_1:
	\Phi(\sigma(j))\in\Omega_2^j\}. \]
Similarly, we define 
\[ \Omega_{2m+1}^j := \{\Phi(\sigma(\ell)):\ell=N_{2m},\ldots,j-1 \}, 
	 \]
\[ N_{2m+1} := \min\left\{ j>N_{2m}: \Phi(\sigma(j))\in \bigcup_{i=0}^{m-1}
	\Omega_{2i+1}^{N_{2i+1}} \cup \Omega_{2m+1}^j \right\}, \]
\[ \Omega_{2m}^j := \{\Phi(\sigma(\ell)):\ell=N_{2m-1},\ldots,j-1 \}, \text{
	and}  
\]
\[ N_{2m} := \min\left\{ j>N_{2m-1}: \Phi(\sigma(j))\in
	\bigcup_{i=1}^{m-1} \Omega_{2m}^{N_{2m}}\cup \Omega_{2m}^j\right\}.
\]

Then there is an $M<q_\alpha$ such that the sets
$\{A_{\sigma(N_i)},\ldots,A_{\sigma(N_{i+1}-1)}\}_{i=0}^{M-1}$ and
$\{A_{\sigma(N_M)},\ldots,A_{\sigma(n)}\}$ separate the $h_j$'s in the sense that 
no $h_j$ is repeated inside any one of these sets, and moreover the union of
``even'' sets contains exactly one copy of each $h_j$ and similarly with the 
``odd'' sets. We can write
\begin{align*}
F_n^\sigma&(k_1\otimes\cdots\otimes k_{p_\alpha}\otimes h_1\otimes
		h_1\otimes\cdots\otimes h_{q_\alpha}\otimes h_{q_\alpha}) \\
	&= \sum_{e_1,\ldots,e_M} \bigg\{
		\langle e_1,[ [[A_{\sigma(1)},A_{\sigma(2)}],\ldots],A_{\sigma(N_1-1)}]\rangle \\
	&\qquad \times
		\left(\prod_{i=1}^{M-1} \langle e_{i+1},
		[ [[e_i,A_{\sigma(N_i)}],\ldots],A_{\sigma(N_{i+1}-1)}]\rangle\right)
		[ [[e_M,A_{\sigma(N_M)}],\ldots],A_{\sigma(n)}]\bigg\}.
\end{align*}
If $M$ is even, then
\begin{multline*} 
\mathcal{A} =
	\langle e_1,[[[A_{\sigma(1)},A_{\sigma(2)}],\ldots],A_{\sigma(N_1-1)}]\rangle \\
	\times \left(\prod_{i=1}^{M/2-1} \langle e_{2i+1},
	[
	[[e_{2i},A_{\sigma(N_{2i})}],\ldots],A_{\sigma(N_{2i+1})-1}]\rangle\right)
	[ [ [e_M, A_{\sigma(N_M)}],\ldots], A_{\sigma(n)}]
\end{multline*}
is a function of $h_1,\ldots,h_{q_\alpha},e_1,\ldots,e_M$ and $k_i$ for $i\in I$ some subset
of $\{1,\ldots,n\}$, and 
\[ \mathcal{B} = 
	\prod_{i=1}^{M/2}\langle e_{2i},
		[[[e_{2i-1},A_{\sigma(N_{2i-1})}],\ldots],A_{\sigma(N_{2i})-1}]\rangle \]
is a function of $h_1,\ldots,h_{q_\alpha},e_1,\ldots,e_M$ and $k_i$ for $i\in
I^c$.  Thus,
\[ 
\|\hat{F}_n^{\sigma,\alpha}\|_2^2
	\le \left(\sum_{k_i:i\in I,h_1,\ldots,h_{q_\alpha},e_1,\ldots,e_M}
			 \|\mathcal{A}\|_{\mathfrak{g}_{CM}}^2\right)
		\left(\sum_{k_i:i\in I^c,h_1,\ldots,h_{q_\alpha},e_1,\ldots,e_M}
			|\mathcal{B}|^2\right) \]
which is finite again by Proposition \ref{p.normest}.
Similarly, if $M$ is odd,
\begin{multline*} 
\mathcal{A} = \langle
		e_1,[[[A_{\sigma(1)},A_{\sigma(2)}],\ldots],A_{\sigma(N_1-1)}]\rangle \\
	\times \prod_{i=1}^{(M-1)/2} \langle e_{2i+1},
		[ [[e_{2i},A_{\sigma(N_{2i})}],\ldots],A_{\sigma(N_{2i+1})-1}]\rangle 
\end{multline*}
and
\[ 
\mathcal{B} = \left(\prod_{i=1}^{(M-1)/2}\langle e_{2i},
		[ [ [e_{2i-1},A_{\sigma(N_{2i-1})}],\ldots],A_{\sigma(N_{2i})-1}]\rangle\right)
		[ [ [e_M, A_{\sigma(N_M)}],\ldots],A_{\sigma(n)}] 
\]
are both functions of $h_1,\ldots,h_{q_\alpha},e_1,\ldots,e_M$ and some $k_i$, and
\[ 
\|\hat{F}_n^{\sigma,\alpha}\|_2^2
	\le \left(\sum_{k_i:i\in I,h_1,\ldots,h_{q_\alpha},e_1,\ldots,e_M}
			 |\mathcal{A}|_{\mathfrak{g}_{CM}}^2\right)
		\left(\sum_{k_i:i\in I^c,h_1,\ldots,h_{q_\alpha},e_1,\ldots,e_M}
			\|\mathcal{B}\|_{\mathfrak{g}_{CM}}^2\right) \]
is again finite by Proposition \ref{p.normest} and this completes the proof.
\end{proof}

Propositions \ref{p.int2} and \ref{p.HS} now allow us to make the following definition.

\begin{defn}
\label{d.BM}
A {\em Brownian motion} on $G$ is the continuous $G$-valued process defined by 
\[g_t = \sum_{n=1}^r \sum_{\sigma\in\mathcal{S}_n}
		\sum_{m=\lceil n/2\rceil}^n \frac{c^\sigma_n }{2^{n-m}}
		\sum_{\alpha\in\mathcal{J}_n^m}
		\int_{\Delta_{p_\alpha}(t)} f_\alpha(s,t) \hat{F}_n^{\sigma,\alpha}
		(dB_{s_1}\otimes\cdots\otimes dB_{s_{p_\alpha}}), \]
where 
\[ c_n^\sigma=(-1)^{e(\sigma)}\bigg/n^2\begin{bmatrix} n-1 \\ e(\sigma) 
	\end{bmatrix}, \]
$\hat{F}_n^{\sigma,\alpha}$ is as defined in equation (\ref{e.Fhat}), 
and $f_\alpha$ is as defined below equation (\ref{e.fprime}).
For $t>0$, let $\nu_t=\mathrm{Law}(g_t)$ be the {\em heat kernel measure at
time $t$}, a probability measure on $G$.
\end{defn}

\begin{prop}[Finite-dimensional approximations]
\label{p.approx}
For $G_\pi$ a finite-dimensional Lie subgroup of
$G_{CM}$, let $\pi$ denote orthogonal projection of
$G_{CM}$ onto
$G_\pi$ and let $g^\pi_t$ be the continuous process on $G_\pi$
defined by
\[ g_t^\pi = \sum_{n=1}^r \sum_{\sigma\in\mathcal{S}_n}
	\sum_{m=\lceil n/2\rceil}^n \frac{c^\sigma_n }{2^{n-m}}
		\sum_{\alpha\in\mathcal{J}_n^m}
	\int_{\Delta_{p_\alpha}(t)} f_\alpha(s,t) \hat{F}_n^{\sigma,\alpha}
	(d\pi B_{s_1}\otimes\cdots\otimes d\pi B_{s_{p_\alpha}}), \]
where the stochastic integrals are defined as in Proposition \ref{p.bad}.
Then $g_t^\pi$ is Brownian motion on $G_\pi$.  In particular, for
$G_\ell=G_{\pi_\ell}$ an increasing sequence of finite-dimensional Lie
subgroups such that the associated orthogonal projections $\pi_\ell$ are
increasing to $I_{\mathfrak{g}_{CM}}$, let
$g^\ell_t=g_t^{\pi_\ell}$.  Then, for all $t<\infty$,
\begin{equation} 
\label{e.b}
\lim_{\ell\rightarrow\infty}\mathbb{E}
	\left\|g^\ell_t-g_t\right\|_\mathfrak{g}^2 = 0. 
\end{equation}
\end{prop}
\begin{proof}
First note that $g_t^\pi$ solves the Stratonovich equation 
$\delta g_t^\pi = L_{g_t^\pi*}\delta \pi B_t$ with $g_0^\pi=\mathbf{e}$, see
\cite{BenArous89,Castell93,Baudoin04} where $\langle \pi B\rangle_t$ is a
standard $\mathfrak{g}_\pi$-valued Brownian motion. 
Thus, $g_t^\pi$ is a $G_\pi$-valued
Brownian motion.

By equation (\ref{e.a}) and its preceding discussion,
\[ g_t^\ell = \sum_{n=1}^r \sum_{\sigma\in\mathcal{S}_n}
	\sum_{m=\lceil n/2\rceil}^n \frac{c^\sigma_n }{2^{n-m}}
		\sum_{\alpha\in\mathcal{J}_n^m}
	\sum_{a=0}^{q_\alpha} b_\alpha^a t^a J_n^\ell(\tilde{f}_\alpha
		\hat{F}_n^{\sigma,\alpha})_t, \]
and thus, to verify (\ref{e.b}), it suffices to show that
\[ \lim_{\ell\rightarrow\infty} \mathbb{E}\|\pi_{\ell}B_t- B_t\|_\mathfrak{g}^2 =
0 \]
and
\[ \lim_{\ell\rightarrow\infty} \mathbb{E} 
	\left\|J_n^\ell(\tilde{f}_\alpha\hat{F}_n^{\sigma,\alpha})_t -
	J_n(\tilde{f}_\alpha\hat{F}_n^{\sigma,\alpha})_t\right\|^2
	 = 0, \]
for all $n\in\{2,\ldots,r\}$, $\sigma\in\mathcal{S}_n$ and
$\alpha\in\mathcal{J}_n$.

So let $\mu_t=\mathrm{Law}(B_t)$.  Then it is known that, if $V$ is a finite-dimensional subspace
of $\mathfrak{g}_{CM}$ and $\pi_V$ is the orthogonal projection from
$\mathfrak{g}_{CM}$ to $V$, then $\pi_V$ admits a $\mu_t$-a.s.~unique
extension to $\mathfrak{g}$.  Moreover, if $V_n$ is an increasing sequence of
finite-dimensional subspaces, then 
\[ \lim_{n\rightarrow\infty} \mathbb{E}\|\pi_{V_n}B_t- B_t\|_\mathfrak{g}^2 = 0;  \]
see for example Section 8.3.3 of \cite{Stroock2011}.

By Proposition \ref{p.HS},
$\hat{F}_n^{\sigma,\alpha}$ is Hilbert-Schmidt, and recall that
$\tilde{f}_\alpha$ is a deterministic polynomial function in $s$.  Thus
$J_n^\ell(\tilde{f}_\alpha \hat{F}_n^{\sigma,\alpha})$ and
$J_n(\tilde{f}_\alpha \hat{F}_n^{\sigma,\alpha})$ are 
$\mathfrak{g}_{CM}$-valued martingales as defined in Proposition \ref{p.bad},
and Proposition \ref{p.bad} gives the desired convergence as well (in
$\mathfrak{g}_{CM}$ and thus in $\mathfrak{g}$).
\end{proof}

\begin{remark}
In fact, for each of the stochastic integrals $J_n(\tilde{f}_\alpha
\hat{F}_n^{\sigma,\alpha})$, it is possible to prove the stronger convergence
that, for all $p\in[1,\infty)$,
\[ \lim_{\ell\rightarrow\infty} \mathbb{E}\left[ \sup_{\tau\le t} 
	\left\|J_n^\ell(\tilde{f}_\alpha\hat{F}_n^{\sigma,\alpha})_\tau -
	J_n(\tilde{f}_\alpha\hat{F}_n^{\sigma,\alpha})_\tau\right\|^p
	\right] = 0, \]
for all $n\in\{2,\ldots,r\}$, $\sigma\in\mathcal{S}_n$ and
$\alpha\in\mathcal{J}_n$. Again, Proposition \ref{p.bad} gives the limit for $p=2$ and thus for
$p\in[1,2]$.  
For $p>2$, Doob's maximal inequality implies it suffices to show that 
\[ \lim_{\ell\rightarrow\infty} \mathbb{E}
	\left\|J_n^\ell(\tilde{f}_\alpha\hat{F}_n^{\sigma,\alpha})_t -
	J_n(\tilde{f}_\alpha\hat{F}_n^{\sigma,\alpha})_t 
	\right\|^p  = 0.  \] 
Since each $J_n^\ell(\tilde{f}_\alpha \hat{F}_n^{\sigma,\alpha})$ and
$J_n(\tilde{f}_\alpha \hat{F}_n^{\sigma,\alpha})$ has chaos expansion
terminating at degree $n$, a theorem of Nelson (see Lemma 2 of
\cite{Nelson73b} and pp. 216-217 of \cite{Nelson73c}) implies that,
for each $j\in\mathbb{N}$, there exists $c_j<\infty$ such that 
\[ \mathbb{E}\left\|J_n^\ell(\tilde{f}_\alpha\hat{F}_n^{\sigma,\alpha})_t -
		J_n(\tilde{f}_\alpha\hat{F}_n^{\sigma,\alpha})_t 
		\right\|^{2j}
	\le c_j \left(\mathbb{E}\left\|
		J_n^\ell(\tilde{f}_\alpha\hat{F}_n^{\sigma,\alpha})_t -
		J_n(\tilde{f}_\alpha\hat{F}_n^{\sigma,\alpha})_t 
		\right\|^2\right)^j. \]
\end{remark}

In a similar way, one may prove the following convergence for the Brownian
motions under right translations by elements of $G_{CM}$.

\begin{prop}
\label{p.rconv}
For any $y\in G_{CM}$, 
\[ \lim_{\ell\rightarrow\infty}
	\mathbb{E}\|g_t^\ell y -g_t y\|^2_{\mathfrak{g}} = 0. \]
where $g_ty$ is the measurable right group action of $y\in G_{CM}$ on $g_t\in
G$, as in Proposition \ref{p.gp}.
\end{prop}

\begin{remark}
	Note that, while the present paper focuses on the case where $\mu$ is non-degenerate and $B$ is Brownian motion on $G$, the above construction and finite-dimensional approximations would all follow with essentially no modification if one considered instead a Gaussian measure $\mu$ whose support was, for example, a subspace $\frak{h}$ of $\frak{g}$ such that $\frak{h}$ generates the span of $\frak{g}$ via the Lie bracket.
\end{remark}

\subsection{Quasi-invariance and log Sobolev}
\label{s.hki}

We are now able to prove Theorem \ref{t.quasi}, which states that the heat kernel measure $\nu_t =
\mathrm{Law}(g_t)$ is
quasi-invariant under left and right translation by elements of $G_{CM}$ and
gives estimates for the Radon-Nikodym derivatives of the ``translated'' measures. Given the results so far, the proof could be given as an application of Theorem 7.3 and Corollary 7.4 of \cite{DriverGordina2009}. However, we provide here a full proof for the reader's convenience.

{\it Proof of Theorem \ref{t.quasi}.}
Fix $t>0$ and $\pi_0$ an orthogonal projection onto a finite-dimensional
subspace $G_0$ of $\mathfrak{g}_{CM}$.  Let $h\in G_0$, and
$\{\pi_n\}_{n=1}^\infty$ be an increasing sequence of projections such that
$G_0\subset \pi_nG_{CM}$ for all $n$ and $\pi_n|_{G_{CM}}\uparrow
I_{G_{CM}}$.  Let $J^{n,r}_t(h,\cdot)$ denote the Radon-Nikodym derivative of
$\nu_t^n\circ R_h^{-1}$ with respect to $\nu_t^n$. Then for each $n$ and for
any $q\in[1,\infty)$, we have the following integrated Harnack inequality
\[ \left(\int_{G_n} \left(J^{n,r}_t(h,g)\right)^q
	d\nu_t^n(g)\right)^{1/q} \le
\exp\left(\frac{(q-1)k}{2(e^{kt}-1)}d_n(e,h)^2\right)
\]
where $k$ is the uniform lower bound on the Ricci curvature as in Proposition
\ref{p.Ric} and $d_n$ is Riemannian distance on $G_n$;
see for example Theorem 1.6 of \cite{DriverGordina2009}.

By Proposition \ref{p.rconv}, we have that for any $f\in C_b(G)$, the class of bounded continuous
functions on $G$
\begin{equation}
\label{e.5.7}
\begin{split}
\int_{G} f(gh) d \nu_{t}(g)
	&= \mathbb{E}[f(g_th)] \\
	&= \lim_{n\rightarrow\infty} \mathbb{E}[f(g_t^nh)]
	= \lim_{n\rightarrow\infty}\int_{G_{n}} (f\circ i_n)(gh) \,
d\nu_t^{n}(g),
\end{split}
\end{equation}
where $ i_n:G_n\rightarrow G$ denotes the inclusion map.  Note that for any $n$
\begin{align*}
\int_{G_n} |(f\circ i_n)(gh)|\,d\nu_t^n(g)
	&= \int_{G_n} J^{n,r}_t(h,g)|(f\circ i_n)(g)|d\nu_t^n(g) \\
	&\le \|f\circ i_n\|_{L^{q'}(G_n,\nu_t^n)} \exp\left(
		\frac{k(q-1)}{2(e^{kt}-1)}d_n(e,h)^2
		\right),
\end{align*}
where $q'$ is the conjugate exponent to $q$.
Allowing $n\rightarrow\infty$ in this last inequality yields
\begin{equation}
\label{e.c}
\int_G |f(gh)|\,d\nu_t(g)
	\le \|f\|_{L^{q'}(G,\nu_t)} \exp\left(
		\frac{k(q-1)}{2(e^{kt}-1)}d(e,h)^2
		\right),
\end{equation}
by equation \eqref{e.5.7} and the fact that
the length of a path in $G_{CM}$
can be approximated by the lengths of paths in the finite-dimensional
projections.  That is, for any
$\pi_0$ and $\varphi\in C^1([0,1],G_{CM})$ with
$\varphi(0)=\mathbf{e}$, there exists an increasing sequence
$\{\pi_n\}_{n=1}^\infty$ of orthogonal projections such that $\pi_0\subset \pi_n$,
$\pi_n|_{\mathfrak{g}_{CM}}\uparrow I_{\mathfrak{g}_{CM}}$, and 
\[ \ell_{CM}(\varphi) = \lim_{n\rightarrow\infty}
	\ell_{G_{\pi_n}}(\pi_n\circ\varphi). \]
To see this, let $\varphi$ be a path in $G_{CM}$.  Then one may show that
\begin{align*}
\ell_{G_{\pi_n}}(\pi_n\circ\varphi)
	&=\int_0^1 \left\|\pi_n\varphi'(s) + \sum_{\ell=1}^{r-1}
		c_\ell \mathrm{ad}_{\pi_n\varphi(s)}^\ell \pi_n\varphi'(s)
		\right\|_{\mathfrak{g}_{CM}}\,ds 
\end{align*}
for appropriate coefficients $c_\ell$; see for example Section 3 of
\cite{Melcher2009}.
Thus, we have proved that \eqref{e.c} holds for $f\in C_b(G)$ and $h\in
\cup_{\pi} G_\pi$.  As this union is dense in $G$ by
Proposition \ref{p.length},
dominated convergence along with the continuity of $d(e,h)$ in $h$ implies
that \eqref{e.c} holds for all $h\in G_{CM}$.

Since the bounded continuous functions are dense in $L^{q'}(G,\nu_t)$ (see for
example Theorem A.1 of \cite{Janson1997}), the inequality in (\ref{e.c}) implies that the
linear functional $\varphi_h:C_b(G)\rightarrow\mathbb{R}$ defined by
\[ \varphi_h(f) = \int_G f(gh)\,d\nu_t(g) \]
has a unique extension to an element, still denoted by
$\varphi_h$, of $L^{q'}(G,\nu_t)^*$ which satisfies the bound
\[ |\varphi_h(f)| \le \|f\|_{L^{q'}(G,\nu_t)}
	\exp\left(
		\frac{k(q-1)}{2(e^{kt}-1)}d(e,h)^2
		\right) \]
for all $f\in L^{q'}(G,\nu_t)$.  Since $L^{q'}(G,\nu_t)^*\cong L^q(G,\nu_t)$, there
then exists a function $J_t^r(h,\cdot)\in
L^q(G,\nu_t)$ such that
\begin{equation}
\label{e.d}
\varphi_h(f) = \int_G f(g)J_t^r(h,g)\,d\nu_t(g),
\end{equation}
for all $f\in L^{q'}(G,\nu_t)$, and
\[ \|J_t^r(h,\cdot)\|_{L^q(G,\nu_t)}
	\le \exp\left(
		\frac{k(q-1)}{2(e^{kt}-1)}d(e,h)^2
		\right). \]

Now restricting (\ref{e.d}) to $f\in C_b(G)$, we may rewrite this equation as
\begin{equation}
\label{e.last}
\int_G f(g)\,d\nu_t(gh^{-1})
	= \int_G f(g) J_t^r(h,g)\,d\nu_t(g).
\end{equation}
Then a monotone class argument (again use Theorem A.1 of
\cite{Janson1997}) shows that (\ref{e.last}) is valid for all
bounded measurable functions $f$ on $G$.  Thus,
$d(\nu_t\circ R_h^{-1})/d\nu_t$ exists and is given by $J_t^r(h,\cdot)$, which is in
$L^q$ for all $q\in(1,\infty)$ and satisfies the desired bound.

A parallel argument gives the
analogous result for $d(\nu_t\circ L_h^{-1})/d\nu_t$.  Alternatively, one
could use the right translation invariance just proved along with the facts that $\nu_t$
inherits invariance under the inversion map $g\mapsto g^{-1}$ from its
finite-dimensional projections and that $d(e,h^{-1})=d(e,h)$.
\hfill$\square$

The following also records the straightforward fact that the heat kernel measure does not charge $G_{CM}$.

\begin{prop}
For all $t>0$, $\nu_t(G_{CM})=0$.
\end{prop}

\begin{proof}
This follows trivially from the fact that $g_t$ is the sum of a Brownian motion $B_t$ on
$\mathfrak{g}$ with a finite sequence of stochastic integrals taking values in
$\mathfrak{g}_{CM}$.
\end{proof}

Thus, $G_{CM}$ maintains its role as a dense subspace of $G$ of
measure 0 with respect to the distribution of the ``group Brownian
motion''.  

\begin{defn}
\label{d.cyl} 
A function $f:G\rightarrow\mathbb{R}$ is said to be a
{\it (smooth) cylinder function} if $f=F\circ\pi$ for some
finite-dimensional projection $\pi$ and
some (smooth) function $F:G_\pi\rightarrow\mathbb{R}$.  Also, $f$ is a 
{\it cylinder polynomial} if $f=F\circ\pi$ 
for $F$ a polynomial function on $G_\pi$.
\end{defn}

\begin{thm}
\label{t.logsob}
Given a cylinder polynomial $f$ on $G$, let
$\nabla f:G\rightarrow\mathfrak{g}_{CM}$ be the gradient of $f$,
the unique element of $\mathfrak{g}_{CM}$ such that
\[ \langle\nabla f(g), h \rangle_{\mathfrak{g}_{CM}} = \tilde{h}f(g)
	:= f'(g)(L_{g*}h_\mathbf{e}), \]
for all $h\in\mathfrak{g}_{CM}$.  
Then for $k$ as in Proposition \ref{p.Ric},
\[ \int_G (f^2\ln f^2)\,d\nu_t -
		\left(\int_G f^2\,d\nu_t \right)\cdot\ln\left(\int_G
		f^2\,d\nu_t\right)
	\le 2\frac{1-e^{-kt}}{k} \int_G \|\nabla f\|_{\mathfrak{g}_{CM}}^2
		\,d\nu_t. \]
\end{thm}

\begin{proof}
Following the method of Bakry and Ledoux applied to $G_P$ (see Theorem 2.9 of 
\cite{Driver96} for the case needed here) shows that
\[ \mathbb{E}\left[\left(f^2\ln f^2\right)\left(g^\pi_t\right)\right] 
		- \mathbb{E}\left[f^2\left(g^\pi_t\right)\right] 
		\ln\mathbb{E}\left[f^2\left(g_t^\pi\right)\right]
	\le 2 \frac{1 - e^{-k_\pi t}}{k_\pi} \mathbb{E}\left\|(\nabla^\pi
		f)\left(g^\pi_t\right)\right\|^2_{\mathfrak{g}_\pi},
\]
for $k_\pi$ as in equation (\ref{e.pah}).  Since the function $x\mapsto
(1-e^{-x})/x$ is decreasing and $k\le k_\pi$ for all finite-dimensional
projections $\pi$, 
this estimate also holds with $k_\pi$ replaced with $k$.  Now applying
Proposition \ref{p.approx} to pass to the limit as $\pi\uparrow I$ gives 
the desired result.
\end{proof}

\bibliographystyle{amsplain}
\bibliography{biblio}

\providecommand{\bysame}{\leavevmode\hbox to3em{\hrulefill}\thinspace}
\providecommand{\MR}{\relax\ifhmode\unskip\space\fi MR }
\providecommand{\MRhref}[2]{%
  \href{http://www.ams.org/mathscinet-getitem?mr=#1}{#2}
}
\providecommand{\href}[2]{#2}
\begin{thebibliography}{10}

\bibitem{AiraultMalliavin2004}
H{\'e}l{\`e}ne Airault and Paul Malliavin, \emph{Backward regularity for some
  infinite dimensional hypoelliptic semi-groups}, Stochastic analysis and
  related topics in {K}yoto, Adv. Stud. Pure Math., vol.~41, Math. Soc. Japan,
  Tokyo, 2004, pp.~1--11. \MR{2083700 (2005h:58062)}

\bibitem{AM2006}
H\'{e}l\`ene Airault and Paul Malliavin, \emph{Quasi-invariance of {B}rownian
  measures on the group of circle homeomorphisms and infinite-dimensional
  {R}iemannian geometry}, J. Funct. Anal. \textbf{241} (2006), no.~1, 99--142.
  \MR{2264248}

\bibitem{ADK2003}
Sergio Albeverio, Alexei Daletskii, and Yuri Kondratiev, \emph{Stochastic
  equations and {D}irichlet operators on infinite product manifolds}, Infin.
  Dimens. Anal. Quantum Probab. Relat. Top. \textbf{6} (2003), no.~3, 455--488.
  \MR{2016321}

\bibitem{AMR00}
Dmitri Alekseevsky, Peter~W. Michor, and Wolfgang~A.F. Ruppert,
  \emph{Extensions of {L}ie algebras}, 2000.

\bibitem{BakhtinMattingly2007}
Yuri Bakhtin and Jonathan~C. Mattingly, \emph{Malliavin calculus for
  infinite-dimensional systems with additive noise}, J. Funct. Anal.
  \textbf{249} (2007), no.~2, 307--353. \MR{2345335 (2008g:60173)}

\bibitem{Baudoin04}
Fabrice Baudoin, \emph{An introduction to the geometry of stochastic flows},
  Imperial College Press, London, 2004. \MR{MR2154760 (2006f:60003)}

\bibitem{BGM2013}
Fabrice Baudoin, Maria Gordina, and Tai Melcher, \emph{Quasi-invariance for
  heat kernel measures on sub-{R}iemannian infinite-dimensional {H}eisenberg
  groups}, Trans. Amer. Math. Soc. \textbf{365} (2013), no.~8, 4313--4350.
  \MR{3055697}

\bibitem{BaudoinTeichmann2005}
Fabrice Baudoin and Josef Teichmann, \emph{Hypoellipticity in infinite
  dimensions and an application in interest rate theory}, Ann. Appl. Probab.
  \textbf{15} (2005), no.~3, 1765--1777. \MR{2152244 (2006g:60080)}

\bibitem{BenArous89}
G{\'e}rard Ben~Arous, \emph{Flots et s\'eries de {T}aylor stochastiques},
  Probab. Theory Related Fields \textbf{81} (1989), no.~1, 29--77. \MR{MR981567
  (90a:60106)}

\bibitem{Bogachev1998}
Vladimir~I. Bogachev, \emph{Gaussian measures}, Mathematical Surveys and
  Monographs, vol.~62, American Mathematical Society, Providence, RI, 1998.
  \MR{MR1642391 (2000a:60004)}

\bibitem{Castell93}
Fabienne Castell, \emph{Asymptotic expansion of stochastic flows}, Probab.
  Theory Related Fields \textbf{96} (1993), no.~2, 225--239. \MR{MR1227033
  (94g:60110)}

\bibitem{CorGrn90}
Lawrence~J. Corwin and Frederick~P. Greenleaf, \emph{Representations of
  nilpotent {L}ie groups and their applications. {P}art {I}}, Cambridge Studies
  in Advanced Mathematics, vol.~18, Cambridge University Press, Cambridge,
  1990, Basic theory and examples. \MR{MR1070979 (92b:22007)}

\bibitem{DalettskiSnaiderman1969}
Ju.~L. Dalecki\u{\i} and Ja.~I. \v{S}na\u{\i}derman, \emph{Diffusion and
  quasiinvariant measures on infinite-dimensional {L}ie groups}, Funkcional.
  Anal. i Prilo\v{z}en. \textbf{3} (1969), no.~2, 88--90. \MR{0248888}

\bibitem{DriverGordina2008}
B.~Driver and M.~Gordina, \emph{Heat kernel analysis on infinite-dimensional
  {H}eisenberg groups}, J. Funct. Anal. \textbf{255} (2008), no.~2, 2395--2461.

\bibitem{DEM2013}
Bruce~K. Driver, Nathaniel Eldredge, and Tai Melcher, \emph{Hypoelliptic heat
  kernels on infinite-dimensional {H}eisenberg groups}, Trans. Amer. Math. Soc.
  \textbf{368} (2016), no.~2, 989--1022. \MR{3430356}

\bibitem{DriverGordina2009}
Bruce~K. Driver and Maria Gordina, \emph{Integrated {H}arnack inequalities on
  {L}ie groups}, J. Differential Geom. \textbf{83} (2009), no.~3, 501--550.
  \MR{MR2581356}

\bibitem{Driver96}
Bruce~K. Driver and Terry Lohrenz, \emph{Logarithmic {S}obolev inequalities for
  pinned loop groups}, J. Funct. Anal. \textbf{140} (1996), no.~2, 381--448.
  \MR{MR1409043 (97h:58176)}

\bibitem{duiskolk}
J.~J. Duistermaat and J.~A.~C. Kolk, \emph{Lie groups}, Universitext,
  Springer-Verlag, Berlin, 2000. \MR{MR1738431 (2001j:22008)}

\bibitem{Elworthy1975}
K.~D. Elworthy, \emph{Measures on infinite-dimensional manifolds}, Functional
  integration and its applications ({P}roc. {I}nternat. {C}onf., {L}ondon,
  1974), 1975, pp.~60--68. \MR{0501086}

\bibitem{FdlP1995}
Denis Feyel and Arnaud de~La~Pradelle, \emph{Brownian processes in infinite
  dimension}, Potential Anal. \textbf{4} (1995), no.~2, 173--183. \MR{1323825}

\bibitem{Gromov1996}
Mikhael Gromov, \emph{Carnot-{C}arath\'eodory spaces seen from within},
  Sub-{R}iemannian geometry, Progr. Math., vol. 144, Birkh\"auser, Basel, 1996,
  pp.~79--323. \MR{1421823 (2000f:53034)}

\bibitem{HairerMattingly2011}
Martin Hairer and Jonathan~C. Mattingly, \emph{A theory of hypoellipticity and
  unique ergodicity for semilinear stochastic {PDE}s}, Electron. J. Probab.
  \textbf{16} (2011), no. 23, 658--738. \MR{2786645 (2012e:60170)}

\bibitem{HuMeyer88-2}
Y.~Z. Hu and P.-A. Meyer, \emph{Chaos de {W}iener et int\'egrale de {F}eynman},
  S\'eminaire de {P}robabilit\'es, {XXII}, Lecture Notes in Math., vol. 1321,
  Springer, Berlin, 1988, pp.~51--71. \MR{MR960508 (89m:60193)}

\bibitem{HuMeyer88-1}
\bysame, \emph{Sur les int\'egrales multiples de {S}tratonovitch}, S\'eminaire
  de {P}robabilit\'es, {XXII}, Lecture Notes in Math., vol. 1321, Springer,
  Berlin, 1988, pp.~72--81. \MR{MR960509 (89k:60070)}

\bibitem{Janson1997}
Svante Janson, \emph{Gaussian {H}ilbert spaces}, Cambridge Tracts in
  Mathematics, vol. 129, Cambridge University Press, Cambridge, 1997.
  \MR{1474726 (99f:60082)}

\bibitem{Kuo1971}
Hui~Hsiung Kuo, \emph{Integration theory on infinite-dimensional manifolds},
  Trans. Amer. Math. Soc. \textbf{159} (1971), 57--78. \MR{295393}

\bibitem{Kuo1972}
\bysame, \emph{Diffusion and {B}rownian motion on infinite-dimensional
  manifolds}, Trans. Amer. Math. Soc. \textbf{169} (1972), 439--459.
  \MR{309206}

\bibitem{Kuo75}
\bysame, \emph{Gaussian measures in {B}anach spaces}, Lecture Notes in
  Mathematics, Vol. 463, Springer-Verlag, Berlin, 1975. \MR{MR0461643 (57
  \#1628)}

\bibitem{Malliavin1990}
Paul Malliavin, \emph{Hypoellipticity in infinite dimensions}, Diffusion
  processes and related problems in analysis, {V}ol.\ {I} ({E}vanston, {IL},
  1989), Progr. Probab., vol.~22, Birkh\"auser Boston, Boston, MA, 1990,
  pp.~17--31. \MR{1110154 (93b:60132)}

\bibitem{Malliavin08}
\bysame, \emph{Invariant or quasi-invariant probability measures for infinite
  dimensional groups: {II}. {U}nitarizing measures or {B}erezinian measures},
  Jpn. J. Math. \textbf{3} (2008), no.~1, 19--47. \MR{MR2390182 (2009b:60170)}

\bibitem{MattinglyPardoux2006}
Jonathan~C. Mattingly and {\'E}tienne Pardoux, \emph{Malliavin calculus for the
  stochastic 2{D} {N}avier-{S}tokes equation}, Comm. Pure Appl. Math.
  \textbf{59} (2006), no.~12, 1742--1790. \MR{2257860 (2007j:60082)}

\bibitem{Melcher2009}
Tai Melcher, \emph{Heat kernel analysis on semi-infinite {L}ie groups}, Journal
  of Functional Analysis \textbf{257} (2009), no.~11, 3552--3592.
  \MR{MR2572261}

\bibitem{Nelson73c}
Edward Nelson, \emph{The free {M}arkoff field}, J. Functional Analysis
  \textbf{12} (1973), 211--227. \MR{MR0343816 (49 \#8556)}

\bibitem{Nelson73b}
\bysame, \emph{Quantum fields and {M}arkoff fields}, Partial differential
  equations ({P}roc. {S}ympos. {P}ure {M}ath., {V}ol. {XXIII}, {U}niv.
  {C}alifornia, {B}erkeley, {C}alif., 1971), Amer. Math. Soc., Providence,
  R.I., 1973, pp.~413--420. \MR{MR0337206 (49 \#1978)}

\bibitem{Pickrell2011}
Doug Pickrell, \emph{Heat kernel measures and critical limits}, Developments
  and trends in infinite-dimensional {L}ie theory, Progr. Math., vol. 288,
  Birkh\"auser Boston, Inc., Boston, MA, 2011, pp.~393--415. \MR{2743770}

\bibitem{Strichartz87}
Robert~S. Strichartz, \emph{The {C}ampbell-{B}aker-{H}ausdorff-{D}ynkin formula
  and solutions of differential equations}, J. Funct. Anal. \textbf{72} (1987),
  no.~2, 320--345. \MR{MR886816 (89b:22011)}

\bibitem{Stroock2011}
Daniel~W. Stroock, \emph{Probability theory}, second ed., Cambridge University
  Press, Cambridge, 2011, An analytic view.

\bibitem{Varadarajan}
V.~S. Varadarajan, \emph{Lie groups, {L}ie algebras, and their
  representations}, Graduate Texts in Mathematics, vol. 102, Springer-Verlag,
  New York, 1984, Reprint of the 1974 edition. \MR{85e:22001}

\end{thebibliography}

\end{document}